\documentclass[12pt]{article}

\usepackage{mathptmx,amsmath,amssymb,amsfonts,amsthm}
\usepackage{graphicx,color,xcolor,pictex,epstopdf,url,algorithm,algorithmic,caption,endnotes}

\usepackage{etoolbox} 
\usepackage{wrapfig} %
\usepackage{mathtools}
\makeatletter
\patchcmd{\@makechapterhead}{\large}{\normalsize}{}{}
\patchcmd{\@makeschapterhead}{\normalsize}{\normalsize}{}{}
\makeatother
\usepackage[textwidth=6.0in,textheight=8.9in,left=0.75in,right=0.75in,top=0.75in,bottom=0.75in]{geometry}
\usepackage[section]{placeins}
\makeatletter
\g@addto@macro\normalsize{\setlength\abovedisplayskip{4pt}}
\g@addto@macro\normalsize{\setlength\belowdisplayskip{4pt}}
\makeatother
\newtheorem{theorem}{Theorem}[section]
\newtheorem{lemma}{Lemma}[section]
\newtheorem{corollary}{Corollary}[section]
\newtheorem{proposition}{Proposition}[section]

\newtheorem{remark}{Remark}[section]
\usepackage{parskip}
\let\oldref\ref
\renewcommand{\ref}[1]{(\oldref{#1})}  
\renewcommand{\eqref}[1]{(\oldref{#1})} 

\newbox\boxaddrone \newbox\boxaddrtwo

\input colordvi
\font\tenrm=cmr10
\font\teni=cmmi10 \skewchar\teni='177
\font\tensy=cmsy10 \skewchar\tensy='60
\font\tenex=cmex10
\font\tenit=cmti10
\font\tensl=cmsl10
\font\tenbf=cmbx10
\font\tentt=cmtt10
\font\ninerm=cmr9
\font\ninei=cmmi9 \skewchar\ninei='177
\font\ninesy=cmsy9 \skewchar\ninesy='60
\font\nineit=cmti9
\font\ninesl=cmsl9
\font\ninebf=cmbx9
\font\ninett=cmtt9
\font\eightrm=cmr8
\font\eighti=cmmi8 \skewchar\eighti='177
\font\eightsy=cmsy8 \skewchar\eightsy='60
\font\eightit=cmti8
\font\eightsl=cmsl8
\font\eightbf=cmbx8
\font\eighttt=cmtt8
\font\sevenrm=cmr7
\font\seveni=cmmi7 \skewchar\seveni='177
\font\sevensy=cmsy7 \skewchar\sevensy='60
\font\sevenbf=cmbx7
\font\sevenit=cmmi7
\font\sevensl=cmmi7
\font\seventt=cmr7
\font\sixrm=cmr6
\font\sixi=cmmi6 \skewchar\sixi='177
\font\sixsy=cmsy6 \skewchar\sixsy='60
\font\sixbf=cmbx6
\font\fiverm=cmr5
\font\fivei=cmmi5 \skewchar\fivei='177
\font\fivesy=cmsy5 \skewchar\fivesy='60
\font\fivebf=cmbx5
\def\tenpoint{\def\rm{\fam0\tenrm}%
        \textfont0=\tenrm \scriptfont0=\sevenrm \scriptscriptfont0=\fiverm
        \textfont1=\teni \scriptfont1=\seveni \scriptscriptfont1=\fivei
        \textfont2=\tensy \scriptfont2=\sevensy \scriptscriptfont2=\fivesy
        \textfont3=\tenex \scriptfont3=\tenex \scriptscriptfont3=\tenex
        \def\it{\fam\itfam\tenit}%
        \textfont\itfam=\tenit
        \def\sl{\fam\slfam\tensl}%
        \textfont\slfam=\tensl
        \def\bf{\fam\bffam\tenbf}%
        \textfont\bffam=\tenbf \scriptfont\bffam=\sevenbf
                \scriptscriptfont\bffam=\fivebf
        \def\tt{\fam\ttfam\tentt}%
        \textfont\ttfam=\tentt
        \normalbaselineskip=12pt%
        \let\sc=\eightrm        
        \setbox\strutbox=\hbox{\vrule height8.5pt depth3.5pt width0pt}%
        \normalbaselines\rm}
\def\ninepoint{\def\rm{\fam0\ninerm}%
        \textfont0=\ninerm \scriptfont0=\sixrm \scriptscriptfont0=\fiverm
        \textfont1=\ninei \scriptfont1=\sixi \scriptscriptfont1=\fivei
        \textfont2=\ninesy \scriptfont2=\sixsy \scriptscriptfont2=\fivesy
        \textfont3=\tenex \scriptfont3=\tenex \scriptscriptfont3=\tenex
        \def\it{\fam\itfam\nineit}%
        \textfont\itfam=\nineit
        \def\sl{\fam\slfam\ninesl}%
        \textfont\slfam=\ninesl
        \def\bf{\fam\bffam\ninebf}%
        \textfont\bffam=\ninebf \scriptfont\bffam=\sixbf
                \scriptscriptfont\bffam=\fivebf
        \def\tt{\fam\ttfam\ninett}%
        \textfont\ttfam=\ninett
        \normalbaselineskip=11pt%
        \let\sc=\sevenrm        
        \setbox\strutbox=\hbox{\vrule height8pt depth3pt width0pt}%
        \normalbaselines\rm}
\def\eightpoint{\def\rm{\fam0\eightrm}%
        \textfont0=\eightrm \scriptfont0=\sixrm \scriptscriptfont0=\fiverm
        \textfont1=\eighti \scriptfont1=\sixi \scriptscriptfont1=\fivei
        \textfont2=\eightsy \scriptfont2=\sixsy \scriptscriptfont2=\fivesy
        \textfont3=\tenex \scriptfont3=\tenex \scriptscriptfont3=\tenex
        \def\it{\fam\itfam\eightit}%
        \textfont\itfam=\eightit
        \def\sl{\fam\slfam\eightsl}%
        \textfont\slfam=\eightsl
        \def\bf{\fam\bffam\eightbf}%
        \textfont\bffam=\eightbf \scriptfont\bffam=\sixbf
                \scriptscriptfont\bffam=\fivebf
        \def\tt{\fam\ttfam\eighttt}%
        \textfont\ttfam=\eighttt
        \normalbaselineskip=9pt%
        \let\sc=\sixrm  
        \setbox\strutbox=\hbox{\vrule height7pt depth2pt width0pt}%
        \normalbaselines\rm}
\def\sevenpoint{\def\rm{\fam0\sevenrm}%
        \textfont0=\sevenrm \scriptfont0=\fiverm \scriptscriptfont0=\fiverm
        \textfont1=\seveni \scriptfont1=\fivei \scriptscriptfont1=\fivei
        \textfont2=\sevensy \scriptfont2=\fivesy \scriptscriptfont2=\fivesy
        \textfont3=\tenex \scriptfont3=\tenex \scriptscriptfont3=\tenex
        \def\it{\fam\itfam\sevenit}%
        \textfont\itfam=\sevenit
        \def\sl{\fam\slfam\sevensl}%
        \textfont\slfam=\sevensl
        \def\bf{\fam\bffam\sevenbf}%
        \textfont\bffam=\sevenbf \scriptfont\bffam=\fivebf
                \scriptscriptfont\bffam=\fivebf
        \def\tt{\fam\ttfam\seventt}%
        \textfont\ttfam=\seventt
        \normalbaselineskip=8pt%
        \let\sc=\fiverm  
        \setbox\strutbox=\hbox{\vrule height6pt depth2pt width0pt}%
        \normalbaselines\rm}

\input pictex
\newbox\figurelegendone
\newbox\figurelegendtwo
\newbox\figureone
\newbox\figuretwo
\newbox\figurethree
\newbox\figurefour
\newbox\figurefive
\newbox\textone
\newbox\texttwo
\newdimen\xfiglen \newdimen\yfiglen  \newdimen\textlen

\begin{document}

\title{Determining the nonlinearity in an acoustic wave equation
}
\author{Barbara Kaltenbacher\footnote{
Department of Mathematics,
Alpen-Adria-Universit\"at Klagenfurt.
barbara.kaltenbacher@aau.at.}
\and
William Rundell\footnote{
Department of Mathematics,
Texas A\&M University,
Texas 77843. 
rundell@math.tamu.edu}
}
\date{\vskip-3ex}
 \maketitle  

\begin{abstract}
\noindent
We consider an undetermined coefficient inverse problem for a nonlinear partial
differential equation describing high intensity ultrasound propagation
as widely used in medical imaging and therapy.
The usual nonlinear term in the standard model using the Westervelt equation in pressure formulation
is of the form $p p_t$. However, this should be considered as a low order
approximation to a more complex physical model where higher order
terms will be required.  Here we assume a more general case where
the form taken is $f(p)\,p_t$ and $f$ is unknown and must be recovered
from data measurements.
Corresponding to the typical measurement setup, the overposed data consists
of time trace observations of the acoustic pressure at a single point or on a one dimensional set $\Sigma$ representing the receiving transducer array at a fixed time.
Additionally to an analysis of well-posedness of the resulting {\sc pde}, we show
injectivity of the linearized forward map from $f$ to the overposed data and use this as motivation for several iterative schemes to recover $f$.
Numerical simulations will also be shown to illustrate the efficiency of the
methods.
\end{abstract}

\leftline{\small \qquad\quad
{\bf Keywords:} Inverse problem, nonlinear acoustics,
reconstruction algorithms}
\smallskip
\leftline{\small \qquad\quad 
\textbf{{\textsc ams}} {\bf classification:} 35R30, 35L15, 35L71.
}

\section{Introduction}
The use of ultrasound is well established in the imaging of human tissue
and the propagation of high intensity ultrasound is modeled by nonlinear wave equations.
Nonlinearity enters the model via the state equation, which is a constitutive relation between acoustic pressure and mass density. This results in the typical situation of the nonlinear effect appearing as a product of a function of the state variable and its time derivative.
A common such model is the Westervelt equation in which the Taylor expansion
of this constitutive relation is truncated to its second degree terms and in this case
a certain ratio of quantities $B/A$ governs the nonlinearity, cf \cite[Chapter 2]{HamiltonBlackstock98}.
\begin{equation}\label{eqn:Westervelt}
u_{tt}-c^2\triangle u -b\triangle u_t = (\kappa u^2)_{tt} = (2\kappa u u_t)_t
\end{equation}
where $c$ is the wave speed, $b$ a damping coefficient and
$\kappa(x)$ is proportional to $B/A$ and may depend on the spatial
variable $x\in {\mathbb {R}}^d$.
Its recovery from overposed data was recently investigated in \cite{nonlinearity_imaging_Westervelt}.

In this paper, our focus will be on using a general nonlinear state equation and on identifying this nonlinearity from indirect measurements.
In section~\oldref{sec:model} we will give more details on the above models
which result from physical laws governing quantities such as the
acoustic particle velocity, the acoustic pressure and the mass density.
Combining these laws while including successive higher order terms,
(the Blackstock scheme) one arrives at a succession of more complex
partial differential equations.
In particular, instead of truncating the Taylor coefficients in the model
we can leave the nonlinearity as a function thus arriving at the equation
\begin{equation}\label{eqn:Westervelt_hyper_fW}
u_{tt} -c^2\triangle u -b\triangle u_t = \bigl(f_W(u) u_t\bigr)_{t}.
\end{equation}
Now the question is whether we can go beyond the finite Taylor series
paradigm of recovering a limited number of lower order coefficients
and seek to recover the constitutive function $f_W$ in
\eqref{eqn:Westervelt_hyper_fW}.

We will work in a finite domain $\Omega\subset {\mathbb {R}}^d$ where
$\Omega$ is a simply connected domain with smooth boundary $\partial\Omega$.
It will turn out to be convenient as in \cite{nonlinearity_imaging_Westervelt}
to rewrite \eqref{eqn:Westervelt_hyper_fW} as a parabolic equation with a
non-local  memory term
\begin{equation*}
u_t -b\triangle u - c^2\int_0^t \triangle u(\tau)\, d\tau = f_W(u)u_t\,.
\end{equation*}
This was the method taken in \cite{nonlinearity_imaging_Westervelt}
and relates the problem to be treated to the problem of the
identification of an unknown nonlinear specific heat coefficient,
see~\cite{PilantRundell:1990}.
The traditional assumption that damping is proportional to velocity
used in \eqref{eqn:Westervelt_hyper_fW} can also be modified to include
so-called fractional damping where the time derivative is replaced by
a subdiffusion operator of order $\alpha$ as, e.g., in \cite{CaiChenFangHolm_survey2018,KaltenbacherRundell:2021a} and the references therein.
In this situation the integral form of the equation is more convenient
as the the standard integral is then replaced by one of Abel fractional type.
We assume that the operator $-\triangle$ on $\Omega$ is equipped with
boundary conditions on $\partial\Omega$ and impose initial conditions.

Typical observations available in such nonlinear acoustic experiments are measurements of the acoustic pressure at an array of transducers or hydrophones.
Thus, there are two obvious types of overposed data: one measuring in the spatial
dimension and the other in the time domain.
In the first case, for some subdomain or curve $\omega$ lying in $\Omega$ and for some
fixed time $T$ we measure
\begin{equation}\label{eqn:final_time_data}
g(x)=u(x,T), \quad x\in \omega\subseteq\Omega
\end{equation}
while in the second for some point $x_0\in\overline{\Omega}$ we measure
\begin{equation}\label{eqn:time_trace_data}
h(t)=u(x_0,t), \quad 0\leq t \leq T.
\end{equation}
As in all problems where the unknown coefficient depends on the state variable
$u$, we must impose range conditions in the above.
Specifically, the values of $u$ occurring in $\Omega\times[0,T]$
must be as subset of those on the measurement domain.
This will naturally impose constraints on the type of initial and boundary
conditions possible.

The outline of the paper is as follows.
We will provide more details on the physical background leading to the
equations in order to justify why the recovery of the term $f_W$ represents
a significant advancement over models relying on fixed Taylor expansions
and (possibly) unknown coefficient terms.
The section following will provide an analysis of the forward operator.
Section~\oldref{sec:reconstruction_schemes} will give the two main reconstruction tools we will use;
a (quasi) Newton scheme to recover $f_W$ using the map $u(x,t;f_W)$ where
this is projected onto the overposed boundary;
fixed point schemes based on an iterative map of Picard type, where the latter will be sped up by Anderson acceleration \cite{Anderson:1965,EvansPollockRebholzXiao:2020,WalkerPeng:2011}.
We will provide an analysis of the forward problem and of well-definedness of the iteration schemes as well as the results of numerical experiments
to demonstrate the effectiveness (and limitations) of these methods for
the problem at hand.

\section{The model}\label{sec:model}
We first derive a generalized version of equation \eqref{eqn:Westervelt}
that results from replacing the second order Taylor expansion underlying the
conventional state equation by a general model function.
Then following the usual elimination steps to obtain a nonlinear second order
wave equation (see, e.g., \cite{HamiltonBlackstock98}), we arrive at a
generalized version of the Westervelt equation.  
To this end, we start with the fundamental quantities of acoustics,
which are the following space and time dependent functions
\begin{itemize}
\item acoustic particle velocity $\vec{v}$;
\item acoustic pressure $p$;
\item mass density $\varrho$;
\end{itemize}
which can be decomposed into their constant mean and a fluctuating part
\[
\vec{v}=\vec{v}_0+\vec{v}_\sim\,, \quad p=p_0+p_\sim\,, \quad \varrho=\varrho_0+\varrho_\sim,
\]
where $\vec{v}_0=0$ in the absence of a flow.

These quantities are interrelated by the following physical balance and material laws:
\begin{itemize}
\item balance of momentum
\begin{equation}\label{NavierStokes}
\varrho \Big(\vec{v}_t + \nabla|\vec{v}|^2\Big) + \nabla p =
\tilde{\mu} \triangle \vec{v}\,;
\end{equation}
\item balance of mass
\begin{equation}\label{conteq}
\nabla\cdot(\varrho\vec{v})=-\varrho_t\,;
\end{equation}
\item state equation relating the acoustic pressure and density fluctuations $p_\sim$ and $\varrho_\sim$:\\
Only in special cases is this known explicitly, an example being perfect gases, with
$\frac{p}{p_0}=\left(\frac{\varrho}{\varrho_0}\right)^\gamma$ where $\gamma$ is the adiabatic index (also known as the ratio of specific heats). The commonly used substitute for obtaining the classical Kuznetsov or Westervelt equations of nonlinear acoustics is  
\begin{equation}\label{stateeq_old}
\varrho_\sim = \frac{p_\sim}{c^2} - \frac{1}{\varrho_0 c^4}\frac{B}{2A} p_\sim^2
- \frac{\tilde{\kappa}}{c^4} {p_\sim}_t\,.
\end{equation}
\end{itemize}
The first two terms in \eqref{stateeq_old} are basically just a polynomial ansatz obtained from a Taylor expansion. We only keep the linear terms, that is, the first term (as this needs to be there to yield a wave type equation) and the last term (as this accounts for attenuation) and  generalize the nonlinearity to be an arbitrary function $\tilde{f}$ in 
\begin{equation}\label{stateeq}
\varrho_\sim = \frac{1}{c^2} p_\sim - \frac{1}{c^4}\tilde{f}(p_\sim)
- \frac{\tilde{\kappa}}{c^4}  {p_\sim}_t\,,
\end{equation}

Subtracting the divergence of \eqref{NavierStokes} from the time derivative of \eqref{conteq} gives
\[
\varrho_{tt}-\triangle p = -\nabla\cdot\Bigl(\varrho_t\vec{v}-\varrho\nabla|\vec{v}|^2+\tilde{\mu} \triangle \vec{v}\Bigr)\,.
\]
Inserting \eqref{stateeq}, and using $\varrho_{0,t}=0$, $\nabla\varrho_0=0$, $p_{0,t}=0$, $\nabla p_0=0$  as well as the approximations 
$\nabla\cdot(\rho_{\sim,t}\vec{v})-\nabla\rho_\sim\cdot\nabla|\vec{v}|^2\approx0$,
$p_{\sim,ttt}\approx c^2\triangle p_{\sim,t}$, 
$\rho_0c^2\nabla\cdot\triangle\vec{v}\approx-\triangle p_{\sim,t}$, 
$\rho_0c^2\triangle |\vec{v}|^2\approx \rho_0|\vec{v}|^2_{tt}\approx \frac{1}{\rho_0^2c^2}(p^2)_{tt}$
that are usually applied in the derivation of the Westervelt equation, setting $b=\kappa+\frac{\mu}{\varrho_0}$
we arrive at
$$
p_{\sim,tt}-c^2\triangle p_\sim -b\triangle p_{\sim,t} = \frac{1}{c^2}\Bigl(\frac{1}{\rho_0^2} p_\sim^2 + \tilde{f}(p_\sim)\Bigr)_{tt}
$$
Now skipping the subscripts $\sim$, and with the abbreviation 
$f(p)= \frac{2}{\rho_0^2} p + \tilde{f}'(p)$, we arrive at the following generalisation of the Westervelt equation
\begin{equation}\label{Westervelt_f_c}
p_{tt}-c^2\triangle p -b\triangle p_t = \frac{1}{c^2}(f(p)p_t)_t\,.
\end{equation}
It is sometimes convenient (e.g., for computational purposes and extending the 
model to fractional operators)
to transform this to a parabolic equation with memory by applying time
integration and using homogeneous initial conditions as well as $\tilde{f}(0)=0$
\begin{equation}\label{Westervelt_f_int_c}
p_t -b\triangle p - c^2\int_0^t \triangle p(\tau)\, d\tau = \frac{1}{c^2} f(p)p_t\,.
\end{equation}
Combining the right hand term with the very first on the left hand side we can rewrite this as
\begin{equation}\label{Westervelt_k_int_c}
(1-\tfrac{1}{c^2}f(p))p_t -b\triangle p - c^2\int_0^t \triangle p(\tau)\, d\tau = 0
\end{equation}
and so the problem is related to the identification of an unknown nonlinear specific heat coefficient, see \cite{PilantRundell:1990}.

Here one might wish to keep track of $c^2$, first of all because of its different order of magnitude as compared to $b$ and $1$, and secondly in order to possibly consider simultaneous identification of $c$ and $f$. In this paper, for simplicity of exposition we just merge it into $f$ by replacing $f\leftarrow \frac{1}{c^2} f$.
Moreover, similarly to \cite{nonlinearity_imaging_Westervelt}, we will add a driving term $r=r(x,t)$ on the right hand side of the {\sc pde} that is supposed to model excitation
and thus arrive at
\begin{equation}\label{Westervelt_f}
p_{tt}-c^2\triangle p -b\triangle p_t = (f(p)p_t)_t + r\,.
\end{equation}
and equivalently 
\begin{equation}\label{Westervelt_k_int}
(1-f(p))p_t -b\triangle p - c^2\int_0^t \triangle p(\tau)\, d\tau = \tilde{r}=\int_0^\cdot r(\tau)\, d\tau + (1-f(p_0))p_1 - b\triangle p_0\,.
\end{equation}

\subsection*{Extension to fractional damping}
In case of fractional damping (more precisely, the Caputo-Wismer model see, e.g., \cite{CaiChenFangHolm_survey2018,KaltenbacherRundell:2021a}, we replace \eqref{Westervelt_f} by
\begin{equation}\label{Westervelt_f_frac}
p_{tt}-c^2\triangle p -b\triangle D^\alpha_t p = \frac{1}{c^2}\Bigl(\frac{1}{\rho_0^2} p^2 + \tilde{f}(p)\Bigr)_{tt} + r
\end{equation}
or after time integration and using homogeneous initial conditions (so it does not matter whether $D^\alpha_t$ is the Djrbashian-Caputo or the Riemann-Liouville fractional derivative) as well as $\tilde{f}(0)=0$ and again using $f(p)= \frac{1}{c^2}(\frac{2}{\rho_0^2} p + \tilde{f}'(p))$
\begin{equation}\label{Westervelt_f_int_frac}
(1-f(p))p_t -b\triangle I_t^{1-\alpha}p - c^2\int_0^t \triangle p(\tau)\, d\tau = \tilde{r}
\end{equation}
which can be rewritten as 
\begin{equation}\label{Westervelt_f_int_frac2}
(1-f(p))p_t -\triangle \int_0^t\Bigl(\tfrac{b}{\Gamma(1-\alpha)}(t-\tau)^{-\alpha} + c^2 \Bigr) p(\tau)\, d\tau = \tilde{r}
\end{equation}

\subsection*{The inverse problem}
The modeling task of reconstructing $f$ in \eqref{Westervelt_f} from overposed data \eqref{eqn:final_time_data} or \eqref{eqn:time_trace_data} thus amounts to inverting the forward map $F=\mbox{tr}_\Sigma\circ G$ which is a composition of the parameter-to-state-map $G:f\mapsto p=p(x,t;f)$ of the initial boundary value problem for \eqref{Westervelt_f} with the trace operator on $\Sigma=\{x_0\}\times (0,T)$ in the time trace and $\Sigma=\omega\times \{T\}$ in the final time case.
Before showing well-definedeness of $F$ in section \oldref{sec:forward} and devising some reconstruction schemes in section \oldref{sec:reconstruction_schemes}, we briefly comment on unique invertibility of its linearisation in the remainder of this section. 

\subsection*{Linearised injectivity}
Linearisation at $f=0$ simplifies the equations to 
\begin{equation}\label{Westervelt_f_PDEinit_lin_0}
\begin{aligned}
&z^0_{tt}-c^2\triangle z^0 -b\triangle z^0_t = (\underline{df}(p^0)p^0_t)_t, \\
&z^0(0)=0, \quad z^0_t(0)=0
\end{aligned}
\end{equation}
where
\begin{equation}\label{Westervelt_f_PDEinit_0}
\begin{aligned}
&p^0_{tt}-c^2\triangle p^0 -b\triangle p^0_t = r, \\ 
&p^0(0)=p_0, \quad p^0_t(0)=p_1.
\end{aligned}
\end{equation}
We can say more about injectivity as well as ill-posedness of the inverse problem linearised at $f=0$ for specially chosen driving and initial functions $r(x,t)=\xi(x)\eta''(t)-b\triangle\xi(x)\eta'(t)-c\triangle\xi(x)\eta(t)$, $p_0(x)=\xi(x)\eta(0)$, $p_1(x)=\xi(x)\eta'(0)$ leading to $p^0(x,t)=\xi(x)\eta(t)$.
Indeed, in  the time trace case, setting $\xi=1$, and choosing $\eta$ such that $\eta(0)=0$ and $\eta'(0)=0$ we can explicitly write the (unique) solution of \eqref{Westervelt_f_PDEinit_lin} as
$z^0(x,t)=\widetilde{\underline{df}}(\eta(t))$ where 
$\widetilde{\underline{df}}(s)=\int_0^s \underline{df}(\tau)\, d\tau$, 
which gives the explicit reconstruction 
$\underline{df}(s)=\frac{h'(\eta^{-1}(s))}{\eta'(\eta^{-1}(s))}$. Ill-posedness therefore results from differentiation of the data and from $\eta'(0)=0$.
\def\ppsi{u}
\section{Analysis of the forward problem}\label{sec:forward}
The question of well-definedeness of the forward operator $F=\mbox{tr}_\Sigma\circ G$ amounts to proving existence of a unique solution $G(f)=p$ to 
\begin{equation}\label{Westervelt_f_PDEinit_forw}
\begin{aligned}
&[(1-f(p))p_t]_t+c^2\mathcal{A} p +b\mathcal{A} p_t = r, \\ 
&p(0)=p_0, \quad p_t(0)=p_1
\end{aligned}
\end{equation}
that is regular enough to admit a trace on $\Sigma=\{x_0\}\times (0,T)$ in the time trace and $\Sigma=\omega\times \{T\}$ in the final time case, respectively.
Here, we denote by $\mathcal{A}$ the negative Laplacian with homogeneous Dirichlet boundary conditions. One could as well incorporate a potentially varying sound speed $c_0(x)$ and other boundary conditions $\mathcal{A}$ as we did in \cite{nonlinearity_imaging_Westervelt},

Setting $\psi(x,t)=\int_0^t p(x,\tau)\, d\tau$ so that $p(x,t)=\psi_t(x,t)$ and integrating \eqref{Westervelt_f_PDEinit_forw} with respect to time we arrive at 
\begin{equation}\label{Westervelt_u}
\begin{aligned}
&(1-f(\psi_t))\psi_{tt}+c^2\mathcal{A} \psi + b\mathcal{A} \psi_t = \tilde{r}, \\ 
&\psi(0)=0, \quad \psi_t(0)=p_0
\end{aligned}
\end{equation}
where 
\begin{equation}\label{eqn:rtil}
\tilde{r}(x,t)=\int_0^t r(x,\tau)\, d\tau + (1-f(p_0))p_1(x) + b\mathcal{A} p_0(x).
\end {equation}
This is a formulation that avoids differentiating $f$.
We will see that this advantage can be preserved in the well-posedness proof at least in the one space dimensional setting.
Physically, $\psi$ corresponds to the acoustic velocity potential with the mass density scaled to unity.
Division of \eqref{Westervelt_u} by $(1-f(\psi_t))$ reveals the fact that we are in fact recovering a nonlinear effective wave speed $\frac{c}{\sqrt{1-f(\psi_t)}}$.

To analyse the nonlinear problem \eqref{Westervelt_u}, we first of all consider the linear one 
\begin{equation}\label{Westervelt_u_lin}
\begin{aligned}
&(1-\sigma)\ppsi_{tt}+c^2\mathcal{A} \ppsi + b\mathcal{A} \ppsi_t -\eta \ppsi_t= \tilde{r}, \\ 
&\ppsi(0)=\ppsi_0, \quad \ppsi_t(0)=\ppsi_1
\end{aligned}
\end{equation}
with given $\sigma(x,t)$ bounded away from $1$, $\sigma(x,t)\leq\overline{\sigma}<1$ so that the coefficient of the second time derivative does not degenerate.
Later on, we will set $\sigma=f(\psi_t)$, $\eta=0$, $\ppsi_0=0$, $\ppsi_1=p_0$ in order to prove well-posedness of \eqref{Westervelt_u} by a fixed point argument, as well as $\sigma=f(\psi_t)$, $\eta=f'(\psi_t)\psi_{tt}$ to investigate its linearization.
Existence of a unique solution to \eqref{Westervelt_u_lin} can be established by the usual Faedo-Galerkin approach of discretisation in space with eigenfunctions of $\mathcal{A}$, deriving energy estimates and taking weak limits. We here only focus on the energy estimates for \eqref{Westervelt_u_lin}.
Some of these can be found scattered over in the cited literature for the special case $f(p)=\kappa p$, however, they will here be tailored to minimize the regularity assumptions of $f$, which is why we provide their derivation in the appendix.
Here and below $C_{X,Y}^\Omega$ denotes the norm of the embedding operator $X(\Omega)\to Y(\Omega)$. 
\begin{lemma} \label{lem:enest_lin}
\begin{itemize}
\item If $\ppsi_1\in H_0^1(\Omega)$, $\tilde{r}\in L^2(0,T;L^2(\Omega))$, $\eta\in L^2(0,T;L^q(\Omega))$ 
with 
\begin{equation}\label{eqn:etasmall}
C(\eta)=
\begin{cases}
(C_{H^1,L^\infty}^\Omega)^2 \|\eta\|_{L^2(0,T;L^2(\Omega))}^2 \ \mbox{ if }d=1\\
(C_{H^1,L^6}^\Omega)^2 \|\eta\|_{L^2(0,T;L^3(\Omega))}^2 \ \mbox{ if }d\in\{2,3\}\,.
\end{cases}
< \frac{b(1-\overline{\sigma})}{4}\,,
\end{equation}
and $\sigma\in L^\infty(\Omega\times(0,T))$ with $\sigma\leq\overline{\sigma}$ a.e., 
then 
\begin{equation}\label{eqn:enest_mult_utt}
\begin{aligned}
&\max\left\{\|\nabla \ppsi_t\|_{L^\infty(0,T;L^2(\Omega))}, \,
\sqrt{(1-\overline{\sigma})C_{lo}(\overline{\sigma})}\|\ppsi_{tt}\|_{L^2(0,T;L^2(\Omega))}
\right\}\\
&\leq \sqrt{bC_{lo}(\overline{\sigma})}\, e^{C(T)}\Bigl(\|\nabla \ppsi_1\|_{L^2(\Omega)}^2 
+ \frac{2}{b(1-\overline{\sigma})}\|\tilde{r}\|_{L^2(0,T;L^2(\Omega))}^2\Bigr)^{1/2}.
\end{aligned}
\end{equation}
with 
\begin{equation}\label{eqn:CT}
C_{lo}(\overline{\sigma})=\frac{2(1-\overline{\sigma})}{b(1-\overline{\sigma})-4C(\eta)}\,, \quad
C(T)= C_{lo}(\overline{\sigma}) \Bigl(c^2+\frac{c^4T}{b}\Bigr) T.
\end{equation}
\item If additionally $\ppsi_1\in H^2(\Omega)$, $\nabla\tilde{r}\in L^2(0,T;L^2(\Omega))$, $\eta\in L^2(0,T;H^1(\Omega))$ with
\begin{equation}\label{eqn:etasmall_mult_Autt}
C(\eta):=(C_{H^2,L^\infty}^\Omega)^2 \|\nabla\eta\|_{L^2(0,T;L^2(\Omega))}^2
+(C_{H^1,L^6}^\Omega)^2 \|\eta\|_{L^2(0,T;L^3(\Omega))}^2
< \frac{b(1-\overline{\sigma})}{4}\,,
\end{equation} 
and $\sigma \in L^\infty(0,T;W^{1,3}(\Omega))$ with 
\begin{equation}\label{eqn:smallnablasigma}
C(\sigma):=C_{H^1\to L^6}^\Omega\|\nabla\sigma\|_{L^\infty(0,T;L^3(\Omega))}< \frac{1-\overline{\sigma}}{2}
\end{equation}
then 
\begin{equation}\label{eqn:enest_mult_Autt}
\begin{aligned}
&\max\left\{\|\mathcal{A} \ppsi_t\|_{L^\infty(0,T;L^2(\Omega))}, \,
\sqrt{(1-\overline{\sigma})C_{hi}(\overline{\sigma})}\|\nabla \ppsi_{tt}\|_{L^2(0,T;L^2(\Omega))}
\right\}\\
&\leq \sqrt{bC_{hi}(\overline{\sigma})}\, e^{C(T)}\Bigl(\|\mathcal{A} \ppsi_1\|_{L^2(\Omega)}^2 
+ \frac{2}{b(1-\overline{\sigma})}\|\nabla\tilde{r}\|_{L^2(0,T;L^2(\Omega))}^2\Bigr)^{1/2}.
\end{aligned}
\end{equation}
with 
\[
C_{hi}(\overline{\sigma})=\frac{2(1-\overline{\sigma}-2C(\sigma))}{b(1-\overline{\sigma})-4C(\eta)}\,, \quad
C(T)= C_{hi}(\overline{\sigma}) \Bigl(c^2+\frac{c^4T}{b}\Bigr) T.
\]
\end{itemize}
\end{lemma}

We first of all consider the spatially 1-d case (which would be less relevant for imaging but is in fact of interest for the problem of recovering the univariate function $f$) since it allows to work with less regular $f$ than the 2- and 3-dimensional setting. 
Indeed, in one space dimension it is enough to use the low level energy estimate \eqref{eqn:enest_mult_utt} and therefore avoid to assume differentiability of $f$. 
Using continuity of the embedding $H^1(\Omega)\to L^\infty(\Omega)$ in one space dimension, the estimate on $\|\nabla \ppsi_t\|_{L^\infty(0,T;L^2(\Omega))}$ in \eqref{eqn:enest_mult_utt} allows us to bound $\|\ppsi_t\|_{L^\infty(\Omega\times(0,T))}$.
By \cite[Lemma 3.3]{temam:2012} we even get $\ppsi_t\in C_w([0,T];C^\beta(\Omega))$ for any $\beta\in[0,\frac12)$, where $C_w$ denotes the space of weakly continuous functions
\begin{equation}\label{eqn:uHoelder}
\begin{aligned}
\|\ppsi_t\|_{L^\infty(\Omega\times(0,T))}
&
\leq \|\ppsi_t\|_{C_w([0,T];C^\beta(\Omega))}\\
&\leq K e^{C(T)}\Bigl(\|\nabla \ppsi_1\|_{L^2(\Omega)}^2+\frac{1}{b(1-\overline{\sigma})}\|\tilde{r}\|_{L^2(0,T;L^2(\Omega))}^2\Bigr)^{1/2}
\end{aligned}
\end{equation}
for some constant $K$ depending on the constants $b$, $c^2$, $\overline{\sigma}$,
$K=K(b,c^2,\overline{\sigma})$ with $K(b,c^2,\overline{\sigma})\to\infty$ as $\overline{\sigma}\nearrow 1$; the estimate likewise holds with $C^\beta$ replaced by $W^{\beta,\infty}$.

Concerning the nonlinear problem, this estimate also tells us that it suffices to assume $f\leq\overline{\sigma}$ and consider $f$ on the interval $[-M,M]$ where $M=M(\overline{\sigma},\underline{\sigma})$ is an upper bound on the right hand side in \eqref{eqn:uHoelder}, that is
\begin{equation}\label{eqn:boundM}
K e^{C(T)}\Bigl(\|\nabla p_0\|_{L^2(\Omega)}+\frac{1}{b(1-\overline{\sigma})}\tilde{R}\Bigr)
\ \leq \, M
\end{equation}
with 
\begin{equation}\label{eqn:Rtil}
\tilde{R}= b\sqrt{T}\|\mathcal{A} p_0\|_{L^2(\Omega)}
+(1+\underline{\sigma})\sqrt{T}\|p_1\|_{L^2(\Omega)}
+ \|r\|_{(H^1(0,T;L^2(\Omega)))^*}
\end{equation}
where $-\underline{\sigma}$ is a lower bound for $f$ so that 
$0\leq 1-\overline{\sigma}\leq 1-f\leq 1+\underline{\sigma}$.
\footnote{For the standard example $f(p)=\kappa p$ we see that on the other hand we need to impose $\kappa M\leq\overline{\sigma}<1$, that is, $M$ must be small enough in order to generate a nonempty set of admissible nonlinearities $f$. In accordance with \eqref{eqn:boundM}, this can always be achieved by making the initial data $p_0$, the driving term $r$, and/or the final time $T$ small enough.}
Fixing $\overline{\sigma}<1$, $\underline{\sigma}\ge0$, we can thus conclude from \eqref{eqn:enest_mult_utt} that for any 
\begin{equation}\label{eqn:calC}
f\in C([-M,M]) \,, \ -\underline{\sigma}\leq f\leq \overline{\sigma}\,,
\end{equation}
the operator $\mathcal{T}:v\to \psi$ where $\psi$ solves \eqref{Westervelt_u_lin} with $\sigma=f(v)$, $\eta=0$, $\tilde{r}$ as in \eqref{eqn:rtil}, is a self-mapping on the closed convex set
\[
\begin{aligned}
\mathcal{M} = \{\phi\in &H^2(0,T;L^2(\Omega))\cap W^{1,\infty}(0,T;H_0^1(\Omega)) \, : \, \|\phi_t\|_{L^\infty(0,T;L^\infty(\Omega))}\leq M
,\\ &\, 
\|\nabla \phi_t\|_{L^\infty(0,T;L^2(\Omega))}\leq m_1
, \ \|\phi_{tt}\|_{L^2(0,T;L^2(\Omega))}\leq m_2 
\},
\end{aligned}
\]
provided 
\begin{equation}\label{eqn:smallness_selfm}
e^{C(T)} \Bigl(\|\nabla p_0\|_{L^2(\Omega)} +\frac{1}{b(1-\overline{\sigma})}\tilde{R}\Bigr)
\leq \max\left\{\frac{m_1}{\sqrt{bC_{lo}(\overline{\sigma})}}, \, m_2\sqrt{\frac{1-\overline{\sigma}}{b}}, \, \frac{M}{K}\right\}\,.
\end{equation}

The set $\mathcal{M}$ is weak* compact in the Banach space $H^2(0,T;L^2(\Omega))\cap W^{1,\infty}(0,T;H_0^1(\Omega))$, which is the dual of a separable space, so that Schauder's Fixed Point Theorem in locally convex topological spaces
\cite{Fan:1952}
provides existence of a fixed point of $\mathcal{T}$ provided $\mathcal{T}$ is weak* continuous.
To prove the latter, we consider an arbitrary sequence $(v_n)_{n\in\mathbb{N}}\subseteq \mathcal{M}$ with weak* limit $v_n\stackrel{*}{\rightharpoonup} v$. Then, due to the above energy estimates, the sequence defined by $\psi_n=\mathcal{T}(v_n)$ lies in $\mathcal{M}$ and therefore has a weakly* convergent subsequence. The limit $\psi^*$ of any weakly* convergent subsequence $(\psi_{n_k})_{k\in\mathbb{N}}$ needs to coincide with the unique weak solution $\psi=\mathcal{T}(v)$ of \eqref{Westervelt_u_lin} with $\sigma=f(v)$, $\eta=0$, $\tilde{r}$ as in \eqref{eqn:rtil}. Indeed, for any $\phi\in C^\infty(0,T;C_0^\infty(\Omega)$ we have 
\begin{equation}\label{eqn:limpsi}
\begin{aligned}
&\int_0^T \int_\Omega \Bigl(\psi^*_{tt}\phi+b\nabla \psi^*_t\cdot\nabla\phi+c^2\nabla \psi^*\cdot\nabla\phi
-f(v_t)\psi^*_{tt}\phi-\tilde{r}\phi\Bigr)\, dx\, dt\\
&=\lim_{k\to\infty}
\int_0^T \int_\Omega \Bigl(\psi_{n_k,tt}\phi+b\nabla \psi_{n_k,t}\cdot\nabla\phi+c^2\nabla \psi_{n_k}\cdot\nabla\phi-f(v_{n_k,t})\psi_{n_k,tt}\phi-\tilde{r}\phi\Bigr)\, dx\, dt=0\,,
\end{aligned}
\end{equation}
where the limit in the nonlinear term follows from Lebesgue's Dominated Convergence Theorem, together with the fact that by continuity of $f$ and $v_{n_k,t}\to v_t$ in $C(\Omega\times(0,T))$ (by the above mentioned continuous embedding) we have pointwise convergence $f(v_{n_k,t}(x,t))\to f(v_t(x,t))$ for all $(x,t)\in \Omega\times(0,T)$ and \\
$\int_0^T\int_\Omega |f(v_{n_k,t})\psi_{n_k,tt}\phi|\, dx\, dt\leq 
\|\psi_{n_k,tt}\|_{L^2(0,T;L^2(\Omega))}\overline{\sigma}\|\phi\|_{L^\infty(\Omega\times(0,T))}<\infty$.

To obtain uniqueness we also prove contractivity of $\mathcal{T}$, which requires more smoothness of $f$, more precisely, Lipschitz continuity. For any $v$, $\tilde{v}\in\mathcal{M}$, the difference $\hat{\psi}=\psi-\tilde{\psi}=\mathcal{T}(v)-\mathcal{T}(\tilde{v})$ solves
\[ 
\begin{aligned}
&(1-f(v_t))\hat{\psi}_{tt}+c^2\mathcal{A}\hat{\psi} + b\mathcal{A} \hat{\psi}_t 
=(f(v_t)-f(\tilde{v}_t)) \tilde{\psi}_{tt} , \\ 
&\hat{\psi}(0)=0, \quad \hat{\psi}_t(0)=0
\end{aligned}
\]
Thus, using \eqref{eqn:enest_mult_utt}, \eqref{eqn:uHoelder} with $\sigma=f(v_t)$, $\eta=0$, $\tilde{r}=(f(v_t)-f(\tilde{v}_t)) \tilde{\psi}_{tt}$, $\ppsi_1=0$, and the fact that $\|\tilde{\psi}_{tt}\|_{L^2(0,T;L^2(\Omega))}\leq m_2$ due to the already shown self-mapping property of $\mathcal{T}$, we obtain
\[
\begin{aligned}
&\max\left\{\frac{1}{\sqrt{bC_{lo}(\overline{\sigma})}}\|\nabla \hat{\psi}_t\|_{L^\infty(0,T;L^2(\Omega))}, \,
\sqrt{\frac{1-\overline{\sigma}}{b}}\|\hat{\psi}_{tt}\|_{L^2(0,T;L^2(\Omega))}, \,
\frac{1}{K}\|\hat{\psi}_t\|_{L^\infty(0,T;W^{\beta,\infty}(\Omega))}\right\}\\
&\leq \sqrt{\frac{2}{b(1-\overline{\sigma})}} \ e^{C(T)} L\,m_2 \|v_t-\tilde{v}_t\|_{L^\infty(0,T;L^\infty(\Omega))},
\end{aligned}
\]
where $L$ is a Lipschitz constant for $f$.
Thus, provided 
\begin{equation}\label{eqn:smallness_contr}
\sqrt{\frac{2K^2}{b(1-\overline{\sigma})}}\ e^{C(T)} L\,m_2 < 1\,,
\end{equation}
the operator $\mathcal{T}$ is a contraction on $\mathcal{M}$ with respect to the norm $|||w|||:=
\|w_t\|_{L^\infty(0,T;W^{\beta,\infty}(\Omega)))}+\|w(0)\|_{L^\infty(\Omega)}$ on $L^\infty(0,T;W^{\beta,\infty}(\Omega))$. The latter is the dual of a separable space and therefore by Alaoglu's Theorem, $\mathcal{M}$ is weak* closed with respect to $L^\infty(0,T;W^{\beta,\infty}(\Omega))$.
Banach's Fixed Point Theorem therefore implies existence of a unique fixed point of $\mathcal{T}$ in $\mathcal{M}$. 

\begin{theorem} \label{th:wellposed_1d}
For $\Omega\subseteq\mathbb{R}^1$, any fixed $\overline{\sigma}<1$, $\underline{\sigma}\geq0$, $M>0$ and for any $f\in C([-M,M])$ such that $-\underline{\sigma}\leq f\leq \overline{\sigma}$ on $[-M,M]$, $p_0\in H_0^1(\Omega)\cap H^2(\Omega)$, $p_1\in L^2(\Omega)$, $r\in (H^1(0,T;L^2(\Omega)))^*$ satisfying \eqref{eqn:smallness_selfm} with $\tilde{R}$ as in \eqref{eqn:Rtil}, there exists a solution $\psi$ of \eqref{Westervelt_u} with \eqref{eqn:rtil}. This solution satisfies the estimate
\[
\begin{aligned}
&\max\left\{\frac{1}{\sqrt{C_{lo}(\overline{\sigma})}}\|\nabla \psi_t\|_{L^\infty(0,T;L^2(\Omega))}, \,
\sqrt{\frac{1-\overline{\sigma}}{b}}\|\psi_{tt}\|_{L^2(0,T;L^2(\Omega))}, \,
\frac{1}{K}\|\psi_t\|_{C_w([0,T];C^\beta(\Omega))}\right\}\\
&\leq e^{C(T)} \Bigl(\|\nabla p_0\|_{L^2(\Omega)}^2+\frac{2}{b(1-\overline{\sigma})}\tilde{R}^2\Bigr)^{1/2}\,.
\end{aligned}
\]
If additionally $f$ is Lipschitz continuous with Lipschitz constant $L$ such that the smallness condition \eqref{eqn:smallness_contr} on $p_0$ and/or $T$ holds, then this solution is unique.
\end{theorem}

From this we conclude that $p=\psi_t\in L^\infty(0,T;H_0^1(\Omega))\cap H^1(0,T;L^2(\Omega))\cap C_w([0,T];C^\beta(\Omega))$ and will therefore infer well-definedness of the forward operator in Corollary \oldref{cor:Fwelldef} below.

In higher space dimensions we need higher Sobolev regularity to guarantee $\psi_t(t)\in L^\infty(\Omega)$ and therewith exclude degeneracy. Hence we will rely on \eqref{eqn:enest_mult_Autt} to establish a self-mapping property of the operator $\mathcal{T}$ on a set that in view of this estimate will be defined by 
\begin{equation}\label{defM_mult_Autt}
\begin{aligned}
\mathcal{M} = \{\phi\in &H^2(0,T;H_0^1(\Omega))\cap W^{1,\infty}(0,T;H^2(\Omega)) \, : \, \|\phi_t\|_{L^\infty(\Omega\times(0,T))}\leq M
,\\ &\, 
\|\mathcal{A} \phi_t\|_{L^\infty(0,T;L^2(\Omega))}\leq m_1
, \ \|\nabla\phi_{tt}\|_{L^2(0,T;L^2(\Omega))}\leq m_2 
\}.
\end{aligned}
\end{equation}
For this purpose we have to bound 
$\|\nabla\sigma\|_{L^\infty(0,T;L^3(\Omega))}$ for $\sigma=f(v_t)$ and $\|\nabla\tilde{r}\|_{L^2(0,T;L^2(\Omega))}$ for $\tilde{r}$ as in \eqref{eqn:rtil}.
To this end, we estimate
\begin{eqnarray}
\nonumber
\|\nabla\sigma\|_{L^\infty(0,T;L^3(\Omega))}
&=& \|f'(v_t)\nabla v_t\|_{L^\infty(0,T;L^3(\Omega))}\leq \|f'\|_{L^\infty(-M,M)} C_{H^1\to L^3}^\Omega\|\mathcal{A} v_t\|_{L^\infty(0,T;L^2(\Omega))}\\
\nonumber
\|\nabla\tilde{r}\|_{L^2(0,T;L^2(\Omega))}
&\leq& b \|\nabla \mathcal{A} p_0\|_{L^2(\Omega)}
+\|\nabla r\|_{(H^1(0,T;L^2(\Omega)))^*} \\
\label{eqn:Rtil_higherdim}
&&+ (1+\underline{\sigma})\sqrt{T}\|\nabla p_1\|_{L^2(\Omega)}
+ \sqrt{T}\|f'\|_{L^\infty(-M,M)} \|p_1 \nabla p_0 \|_{L^2(\Omega)} =:\tilde{R}
\end{eqnarray}
and assume
\begin{equation}\label{eqn:smallness_nablasigma} 
L C_{H^1\to L^3}^\Omega C_{H^1\to L^6}^\Omega m_1\leq \frac{1-\overline{\sigma}}{4}
\end{equation}
to satisfy \eqref{eqn:smallnablasigma} (note that with $\eta=0$, \eqref{eqn:etasmall_mult_Autt} is anyway satisfied)
\begin{equation}\label{eqn:smallness_selfm_higherdim} 
e^{C(T)} \Bigl(\|\mathcal{A} p_0\|_{L^2(\Omega)} +\frac{1}{b(1-\overline{\sigma})}\tilde{R}\Bigr)
\leq \max\left\{\frac{m_1}{\sqrt{2}}, \, m_2\sqrt{\frac{1-\overline{\sigma}}{2b}}, \, \frac{M}{K}\right\}
\end{equation}
with $\tilde{R}$ as in \eqref{eqn:Rtil_higherdim}.

Since \eqref{eqn:enest_mult_utt} remains valid also in higher space dimensions for $\eta=0$,  contractivity with respect to the $|||\cdot|||$ norm for Lipschitz continuous $f$ with Lipschitz constant satisfying \eqref{eqn:smallness_contr} directly carries over to the spatially higher dimensional setting.

\begin{theorem} \label{th:wellposed}
For $\Omega\subseteq\mathbb{R}^d$, $d\in\{1,2,3\}$, $\partial\Omega\in C^{2,\alpha}$, any fixed $\overline{\sigma}<1$, $\underline{\sigma}\geq0$, $M>0$ and for any $f\in W^{1,\infty}(-M,M)=C^{0,1}([-M,M])$ such that $-\underline{\sigma}\leq f\leq \overline{\sigma}$ on $[-M,M]$, $p_0\in H_0^1(\Omega)\cap H^3(\Omega)$, $p_1\in H^1(\Omega)$, $\nabla r\in (H^1(0,T;L^2(\Omega)))^*$ satisfying \eqref{eqn:smallness_contr}, \eqref{eqn:smallness_nablasigma}, \eqref{eqn:smallness_selfm_higherdim} with $\tilde{R}$ as in \eqref{eqn:Rtil_higherdim}, there exists a unique solution $\psi$ of \eqref{Westervelt_u} with \eqref{eqn:rtil}. This solution satisfies the estimate
\[
\begin{aligned}
&\max\left\{\frac{1}{\sqrt{C_{hi}(\overline{\sigma})}}\|\mathcal{A} \psi_t\|_{L^\infty(0,T;L^2(\Omega))}, \,
\sqrt{\frac{1-\overline{\sigma}}{b}}\|\nabla\psi_{tt}\|_{L^2(0,T;L^2(\Omega))}, \,
\frac{1}{K}\|\psi_t\|_{L^\infty(0,T;W^{\beta,\infty}(\Omega))}\right\}\\
&\leq e^{C(T)} \Bigl(\|\mathcal{A} p_0\|_{L^2(\Omega)}^2+\frac{2}{b(1-\overline{\sigma})}\tilde{R}^2\Bigr)^{1/2}\,.
\end{aligned}
\]
\end{theorem}

With $p=\psi_t$ we thus get 
\begin{equation}\label{eqn:defU}
p\in U:=L^\infty(0,T;H_0^1(\Omega)\cap H^2(\Omega))\cap H^1(0,T;H_0^1(\Omega))\cap C_w([0,T];C^\beta(\Omega))\,.
\end{equation}

\begin{corollary}\label{cor:Fwelldef}
Let $\overline{\sigma}<1$, $\underline{\sigma}\geq0$, $M>0$, $L>0$, $\mathcal{D}=\{f\in W^{1,\infty}(-M,M)\, : \, -\underline{\sigma}\leq f\leq \overline{\sigma}, \ |f'|\leq L \mbox{ a.e. }\}$, and either 
\begin{itemize}
\item[(a)] $\Omega\subseteq\mathbb{R}^1$, 
$p_0\in H_0^1(\Omega)\cap H^2(\Omega)$, $p_1\in L^2(\Omega)$, $r\in (H^1(0,T;L^2(\Omega)))^*$ satisfying  \eqref{eqn:smallness_selfm}, \eqref{eqn:smallness_contr} with $\tilde{R}$ as in \eqref{eqn:Rtil} or
\item[(b)] $\Omega\subseteq\mathbb{R}^d$, $d\in\{2,3\}$
$p_0\in H_0^1(\Omega)\cap H^3(\Omega)$, $p_1\in H_0^1(\Omega)$, $\nabla r\in (H^1(0,T;L^2(\Omega)))^*$ satisfying \eqref{eqn:smallness_nablasigma}, \eqref{eqn:smallness_selfm_higherdim}, \eqref{eqn:smallness_contr} with $\tilde{R}$ as in \eqref{eqn:Rtil_higherdim}.
\end{itemize}
Then the forward operator $F:\mathcal{D}\to Y$ is well-defined, where 
\begin{equation}\label{eqn:defY}
Y=L^{\rm p}(0,T)\mbox{ for time trace observations and } 
Y=L^{\rm p}(\omega)\mbox{ for final time observations,} 
\end{equation}
respectively, for any ${\rm p}\in[1,\infty]$.
\end{corollary}

The following continuity result on $F$ is useful for, e.g., proving that Tikhonov regularisation is well-defined, but we will also use it to establish weak sequential convergence of our iterative reconstruction schemes, cf. Remarks \oldref{rem:convfp} and \oldref{rem:convNewton} below.
 
\begin{proposition}\label{prop:Fcont}
Under the conditions of Corollary \oldref{cor:Fwelldef}, the operator $F:\mathcal{D}\subseteq W^{1,\infty}(-M,M)\to Y$ is weakly(*) continuous, that is, for any sequence $(f_n)_{n\in\mathbb{N}}\subseteq\mathcal{D}$ converging weakly* in $W^{1,\infty}(-M,M)$ to $f$, we have $G(f_n)\stackrel{*}{\rightharpoonup} G(f)$ in $U$ and $F(f_n)\stackrel{*}{\rightharpoonup} F(f)$ in $Y$
(see \eqref{eqn:defU}, \eqref{eqn:defY} for the defitions of $U$ and $Y$).
\end{proposition}
{\it Proof.}
Define $\tilde{U}:=W^{1,\infty}(0,T;H_0^1(\Omega)\cap H^2(\Omega))\cap H^2(0,T;H_0^1(\Omega))\cap \{\psi\, : \, \psi_t\in C_w([0,T];C^\beta(\Omega))\}$.
We proceed by proving that $\tilde{G}:W^{1,\infty}(-M,M)\to\tilde{U}$
is weakly(*) continuous.
The result then follows from linearity and boundedness of $\partial_t:\tilde{U}\to U$ and $\mbox{tr}_\Sigma:U\to Y$.
Let $(f_n)_{n\in\mathbb{N}}\subseteq\mathcal{D}$ be an arbitrary sequence with $f_n\stackrel{*}{\rightharpoonup}f$ in $W^{1,\infty}(-M,M)$. Then $f\in\mathcal{D}$ and by Theorem \oldref{th:wellposed} $(\psi_n)_{n\in\mathbb{N}}:=(\mathcal{G}(f_n))_{n\in\mathbb{N}}$ is bounded in $\tilde{U}$ (more precisely, contained in $\mathcal{M}$). Hence there exists a subsequence $(\psi_{n_k}, f_{n_k})_{k\in\mathcal{N}}$ and an element $(\psi^*,f^*)\in \tilde{U}\times W^{1,\infty}(-M,M)$ such that $(\psi_{n_k},f_{n_k})\stackrel{*}{\rightharpoonup}\psi^*$ in $\tilde{U}\times W^{1,\infty}(-M,M)$, $(\psi_{n_k,t},f_{n_k})\to (\psi^*_t,f^*)$ in $L^2(0,T;L^2(\Omega))\times C([-M,M])$ (by the Arzel\'{a}-Ascoli Theorem), where by uniqueness of limits, $f^*=f$.
Thus, for any $\phi\in C_0^\infty(0,T;C_0^\infty(\Omega)$, similarly to \eqref{eqn:limpsi} we have 
\begin{equation}\label{eqn:limpsi_cont}
\begin{aligned}
&\int_0^T \int_\Omega \Bigl(\psi^*_{tt}\phi+b\nabla \psi^*_t\cdot\nabla\phi+c^2\nabla \psi^*\cdot\nabla\phi
-f(\psi^*_t)\psi^*_{tt}\phi-\tilde{r}\phi\Bigr)\, dx\, dt\\
&=\lim_{k\to\infty}
\int_0^T \int_\Omega \Bigl(\psi_{n_k,tt}\phi+b\nabla \psi_{n_k,t}\cdot\nabla\phi+c^2\nabla \psi_{n_k}\cdot\nabla\phi-f_{n_k}(\psi_{n_k,t})\psi_{n_k,tt}\phi-\tilde{r}\phi\Bigr)\, dx\, dt=0\,,
\end{aligned}
\end{equation}
where for the convergence of the nonlinear term we argue as follows. 
We decompose
\[
\begin{aligned}
& 
f(\psi^*_t)\psi^*_{tt}-f_{n_k}(\psi_{n_k,t})\psi_{n_k,tt}
\\& = 
(f(\psi^*_t)-f(\psi_{n_k,t}))\psi^*_{tt}
+(f(\psi_{n_k,t})-f_{n_k}(\psi_{n_k,t}))\psi^*_{tt}
+f_{n_k}(\psi_{n_k,t})(\psi^*_{tt}-\psi_{n_k,tt})
\end{aligned}
\]
and consider the limit in each of the terms (integrated against $\phi$ over $\Omega\times(0,T)$) separately.
For the first term, we can conclude $\int_0^T \int_\Omega (f(\psi^*_t)-f(\psi_{n_k,t}))\psi^*_{tt} \phi \, dx\, dt\to0$ from Lebesgue's Dominated Convergence Theorem and continuity of $f$, similarly to the proof of \eqref{eqn:limpsi} above.
The second term can be estimated by 
$|\int_0^T \int_\Omega (f(\psi_{n_k,t})-f_{n_k}(\psi_{n_k,t}))\psi^*_{tt} \phi \, dx\, dt|\leq
\|f-f_{n_k}\|_{C(-M,M)} \|\psi^*_{tt} \phi\|_{L^1(0,T;L^1(\Omega)}\to0$.
Finally, we have 
\[
\begin{aligned}
&\left|\int_0^T \int_\Omega
f_{n_k}(\psi_{n_k,t})(\psi^*_{tt}-\psi_{n_k,tt})
\phi \, dx\, dt\right|\\
&\qquad=
\left|-\int_0^T \int_\Omega
\Bigl(f_{n_k}'(\psi_{n_k,t})\psi_{n_k,tt}\phi + f_{n_k}(\psi_{n_k,t})\phi_t\Bigr) (\psi^*_{t}-\psi_{n_k,t})
\, dx\, dt\right|\\
&\qquad\leq \Bigl(L m_2 \|\phi\|_{L^\infty(0,T;L^\infty(\Omega))}+\max\{-\underline\sigma,\overline\sigma\}\|\phi_t\|_{L^2(0,T;L^2(\Omega))}\Bigr) \|\psi^*_{t}-\psi_{n_k,t}\|_{{L^2(0,T;L^2(\Omega))}}\to0\,.
\end{aligned}
\]
From \eqref{eqn:limpsi_cont}, together with a subsequence-subsequence argument and the uniqueness part of Theorem \oldref{th:wellposed} we conclude convergence of the whole sequence $(\psi_n)_{n\in\mathbb{N}}$ to $G(f)$.
\hfill $\diamondsuit$\break

\begin{remark}
The fractional attenuation case \eqref{Westervelt_f_frac} can be tackled similarly in principle.
However, the resulting reformulation via the acoustic velocity potential (that allowed us to minimize smoothness assumptions on $f$ above) does not work any more for $\alpha<1$. More precisely, the natural approach of multiplying the linearized equation  
\[
\begin{aligned}
&(1-\sigma)\ppsi_{tt}+c^2\mathcal{A} \ppsi + b\mathcal{A} D^\alpha_t\ppsi -\eta \ppsi_t= \tilde{r}, \\ 
&\ppsi(0)=\ppsi_0, \quad \ppsi_t(0)=\ppsi_1
\end{aligned}
\]
with $D_t^{1+\alpha}\ppsi(\tau)=D_t(D_t^\alpha\ppsi)(\tau))-\frac{\tau^{-\alpha}}{\Gamma(1-\alpha)}\ppsi_1$ and integrating over $(0,t)$ yields the energy identity 
\[
\begin{aligned}
&\int_0^t \langle \ppsi_{tt}(\tau),(I^{1-\alpha}\ppsi_{tt})(\tau)\rangle_{L^2(\Omega)}\, d\tau
+\frac{b}{2}\|\nabla (D_t^\alpha\ppsi)(t)\|_{L^2(\Omega)}\\
&=\frac{b}{2}\|\nabla (D_t^\alpha\ppsi)(0)\|_{L^2(\Omega)}
+\int_0^t \langle \sigma\ppsi_{tt}(\tau),D_t^{1+\alpha}\ppsi(\tau)\rangle_{L^2(\Omega)}\, d\tau
+ b\int_0^t \tfrac{\tau^{-\alpha}}{\Gamma(1-\alpha)}\langle \ppsi_t(0),D_t^{1+\alpha}\ppsi(\tau)\rangle_{L^2(\Omega)}\, d\tau\\
&\quad+c^2 \int_0^t \tfrac{\tau^{-\alpha}}{\Gamma(1-\alpha)}\langle\nabla\ppsi_t(0),\nabla \ppsi(\tau)\rangle_{L^2(\Omega)}\, d\tau 
+c^2 \int_0^t \langle\nabla\ppsi_t(\tau),\nabla D_t^\alpha\ppsi(\tau)\rangle_{L^2(\Omega)}\, d\tau \\
&\quad-c^2 \langle\nabla\ppsi(t),\nabla D_t^\alpha\ppsi(t)\rangle_{L^2(\Omega)}
+c^2 \langle\nabla\ppsi(0),\nabla D_t^\alpha\ppsi(0)\rangle_{L^2(\Omega)}
+\int_0^t \langle \tilde{r}(\tau),D_t^{1+\alpha}\ppsi(\tau)\rangle_{L^2(\Omega)}\, d\tau\,.
\end{aligned}
\]
Here the first term on the left hand can be estimated from below by means of coercivity of the Abel integral operator 
\cite{VogNedSau:16}
\[
\int_0^t \langle \ppsi_{tt}(\tau),(I^{1-\alpha}\ppsi_{tt})(\tau)\rangle_{L^2(\Omega)}\, d\tau
\geq \cos((1-\alpha)\pi/2)\|\ppsi_{tt}\|_{H^{-(1-\alpha)/2}(0,t;L^2(\Omega))}\,.
\]
the second term on the right hand side can be bounded from above by 
\[
\begin{aligned}
\int_0^t \langle \sigma\ppsi_{tt}(\tau),D_t^{1+\alpha}\ppsi(\tau)\rangle_{L^2(\Omega)}\, d\tau
\leq \|\ppsi_{tt}\|_{H^{-(1-\alpha)/2}(0,t;L^2(\Omega))}
\|\sigma D_t^{1+\alpha}\ppsi\|_{H^{(1-\alpha)/2}(0,t;L^2(\Omega))}
\end{aligned}
\]
where the second factor can be estimated by means of the Kato-Ponce inequality 
\begin{equation}\label{prodruleest}
\|fg\|_{W^{\rho,r}(0,T)}\lesssim \|f \|_{W^{\rho,p_1}(0,T)} \|g\|_{L^{q_1}(0,T)}
+ \| f \|_{L^{p_2}(0,T)} \|g\|_{W^{\rho,q_2}(0,T)}
\end{equation}
for $0\leq\rho\leq\overline{\rho}<1$, 
$1<r<\infty$, $p_1$, $p_2$, $q_1$, $q_2\in(1,\infty]$, with $\frac{1}{r}=\frac{1}{p_i}+\frac{1}{q_i}$, $i=1,2$; see, e.g., \cite{GrafakosOh2014} as follows
\[
\begin{aligned}
\|\sigma D_t^{1+\alpha}\ppsi\|_{H^{(1-\alpha)/2}(0,t;L^2(\Omega))}
\lesssim &
\|\sigma \|_{H^{(1-\alpha)/2}(0,t;L^\infty(\Omega))}
\|D_t^{1+\alpha}\ppsi(\tau)\|_{L^2(0,t;L^2(\Omega))}\\
&+\|\sigma \|_{L^2(0,t;L^2(\Omega))}
\|D_t^{1+\alpha}\ppsi\|_{H^{(1-\alpha)/2}(0,t;L^2(\Omega))}
\end{aligned}
\]
where 
\[
\|D_t^{1+\alpha}\ppsi(\tau)\|_{L^2(0,t;L^2(\Omega))}
\lesssim \|\ppsi_{tt}\|_{H^{-(1-\alpha)}(0,t;L^2(\Omega))}
\]
and 
\[
\|D_t^{1+\alpha}\ppsi\|_{H^{(1-\alpha)/2}(0,t;L^2(\Omega))}
\lesssim \|\ppsi_{tt}\|_{H^{-(1-\alpha)/2}(0,t;L^2(\Omega))}\,.
\]
It can thus be absorbed into the first left hand side term, provided $\sigma$ is small enough.
However, the fifth term on the right hand side yields a nonnegative term
\[
c^2 \int_0^t \langle\nabla\ppsi_t(\tau),\nabla D_t^\alpha\ppsi(\tau)\rangle_{L^2(\Omega)}\, d\tau
=c^2 \int_0^t \langle\nabla\ppsi_t(\tau),\nabla (I^{1-\alpha}\ppsi_t)(\tau)\rangle_{L^2(\Omega)}\, d\tau
\]
on the right hand side of the energy identity that cannot be dominated by any of the nonnegative left hand side energy contributions. This difficulty has already been observed previously, see \cite[Remark 1]{fracJMGT}. 
Therefore, one would have to work with the pressure formulation \eqref{Westervelt_f_PDEinit_forw} 
analogously to \cite[Section 3.1]{KaltenbacherRundell:2021d}. However, this obviously requires higher differentiability of $f$.
\end{remark}

\section{Reconstruction Schemes}\label{sec:reconstruction_schemes}
Now we turn to the inverse problem of identifying $f$ in 
\begin{equation}\label{Westervelt_f_PDEinit}
\begin{aligned}
&[(1-f(p))p_t]_t+c^2\mathcal{A} p +b\mathcal{A} p_t = r, \\ 
&p(0)=p_0, \quad p_t(0)=p_1
\end{aligned}
\end{equation}
from time trace 
\begin{equation}\label{eqn:titr}
h(t)=p(x_0,t), \quad t\in(0,T)
\end{equation}
or final time
\begin{equation}\label{eqn:fiti}
g(x)=p(x,T), \quad x\in\omega\subseteq\Omega
\end{equation}
observations.

\subsection{Fixed point iterations}\label{subsec:Pic}
Fixed point formulations of the inverse problem can be obtained by projecting the {\sc pde} on the observation manifold and inserting the available measurement data where possible.
From this we obtain an iterative reconstruction scheme by applying Picard iteration.

In the case of time trace data \eqref{eqn:titr} using the fact that since $p=p(x,t;f)$ solves \eqref{Westervelt_k_int} we have the identity
\[
\begin{aligned}
(1-f(h(t)))h'(t)&= -\mathcal{A} \Bigl( b p(x_0,t;f) + c^2\int_0^t p(x_0,\tau,f)\, d\tau \Big) \\
&= (1-f(p(x_0,t;f)))p_t(x_0,t;f)
\end{aligned}
\]
and we obtain the fixed point scheme
\begin{equation}\label{iteration_titr}
f_{k+1}(h(t))= 1-\frac{p_t(x_0,t;f_k)}{h'(t)}(1-f_k(p(x_0,t;f_k))))\,.
\end{equation}

For final time data \eqref{eqn:fiti} using the identity
\[
\begin{aligned}
\bigl(1-f(g(x))\bigr)p_t(x,T;f)&=-\mathcal{A} \Bigl( b g(x) + c^2\int_0^t p(x,\tau;f)\, d\tau \Big)\\
&= -b\mathcal{A} (g(x)- p(x,T;f)) + (1-f(p(x,T;f))p_t(x,T;f)
\end{aligned}
\]
we get the fixed point scheme 
\begin{equation}\label{iteration_fiti}
f_{k+1}(g(x))= f_k(p(x,T;f_k)) +\frac{b}{p_t(x,T;f_k)}\mathcal{A} (g(x)- p(x,T;f_k))\,.
\end{equation}

With fractional damping \eqref{Westervelt_f_int_frac2}, 
the time trace data iteration scheme \eqref{iteration_titr} remains exactly the same as in the strong damping case $\alpha=1$, just the {\sc pde} to be solved in between is modified.
For final time data, the possibility of inserting observations becomes very limited, since except for the very first term in \eqref{Westervelt_f_int_frac2}, $p$ only appears under an integral over time; so we do not pursue the final time data case in the context of \eqref{Westervelt_f_int_frac2}.

\subsection{Newton type schemes}\label{subsec:Newt}
Recall that the forward operator $F=\mbox{tr}_\Sigma\circ G$ is a composition of the parameter-to-state map $G:f\mapsto p=p(x,t;f)$ where $p$ solves \eqref{Westervelt_f_PDEinit}  with the trace operator on $\Sigma=\{x_0\}\times (0,T)$ in the time trace and $\Sigma=\omega\times \{T\}$ in the final time case.

The linearisation $z=G'(f)\underline{df}$ solves
\begin{equation}\label{Westervelt_f_PDEinit_lin}
\begin{aligned}
&[(1-f(p))z_t]_t+c^2\mathcal{A} z +b\mathcal{A} z_t - [f(p)'z p_t]_t= (\underline{df}(p)p_t)_t, \\
&z(0)=0, \quad z_t(0)=0\,.
\end{aligned}
\end{equation}

Newton's method is therefore defined by the equations
\begin{equation}\label{eqn:Newt}
\begin{aligned}
&z(x_0,t;f_k)= h(t)- p(x_0,t;f_k), \quad t\in(0,T)\mbox{ in the time trace case} \\
&z(x,T;f_k)= g(x)- p(x,T;f_k), \quad x\in \omega\mbox{ in the final time case} 
\end{aligned}
\end{equation}
where $p(\cdot;f_k)$ solves \eqref{Westervelt_f_PDEinit} with $f=f_k$ and $z(\cdot;f_k)$ solves \eqref{Westervelt_f_PDEinit_lin} with $f=f_k$, $\underline{df}=f_{k+1}-f_k$.

We will also consider a frozen version of Newton's method, where we linearize at a fixed initial guess, that is, rely on $z=G'(f_0)(f_{k+1}-f_k)$ rather than $z=G'(f_k)(f_{k+1}-f_k)$.
Therefore, in the frozen version of \eqref{eqn:Newt}, $p(\cdot;f_k)$ still solves \eqref{Westervelt_f_PDEinit} with $f=f_k$ but $z(\cdot;f_k)$ now solves \eqref{Westervelt_f_PDEinit_lin} with $f=f_0$, $\underline{df}=f_{k+1}-f_k$.

These Newton type methods extend to the case of fractional damping \eqref{Westervelt_f_frac} in a straightforward manner by replacing $b\mathcal{A} z_t$ by $b\mathcal{A} D^\alpha_t z$ in \eqref{Westervelt_f_PDEinit_lin}.

\subsection{Andersen acceleration}\label{subsec:And}
The fact that Picard iteration, namely finding a fixed point of the map
$\mathcal{T}$, using $x_{n+1} = \mathcal{T}(x_n)$,
can be extremely slow is legendary.
The usual approach is to obtain a sufficiently small bound on $\mathcal{T}'$
that the contraction mapping can be used. Even here we are faced with
linear convergence.
There are a variety of methods that have been used to speed up
the convergence of the sequence of iterates and one that is commonly
taken is Anderson acceleration. 
See, for example, \cite{Anderson:1965,WalkerPeng:2011}.

The algorithm with depth $m$ and {\it damping factors\/} $\{\beta_k\}$
is as follows.

\begin{enumerate}
\item
Initial approximation $x_0$ then compute $\tilde x_1 = \mathcal{T}(x_0)$
and set $x_1 = \tilde x_1$.
\item
For $k=1,\,2,\,\ldots $
\begin{itemize}
\item set $m_k = \hbox{\rm min}\{k,m\}$.
\item compute $\tilde x_{k+1} = \mathcal{T}(x_k)$
\item Solve the minimization problem for $\{a_j\}$
(where $a_j$ depends on $k$)
$$ \min\Bigl\{
\Bigl\|\sum_{j=k-m_k}^k a_j(\tilde x_{j+1} - x_j)\Bigr\|
\, : \, \sum_{j=k-m_k}^k a_j=1 \Bigr\}
$$
\item Then set 
$$ x_{k+1} = \sum_{j=k-m_k}^k a_j\bigl((1-\beta_k)x_j +\beta_k \tilde x_{j+1}\bigr)
$$
\end{itemize}
\end{enumerate}

{\bf Some notes/remarks:}
The key is the minimization step. 
Since $m$ is typically small ($3\leq m\leq 5$)
the simplest approach is to use least squares in the sense of 
$w(k):= \tilde x_k - x_{k-1}$ and
\begin{equation}\label{LS_Anderson}
\small
A = \left[\begin{array}{c}w(k)\\ \ldots\\ w(k-m_k)\end{array}\right]
\qquad C = [1\;1\,\ldots \;1]\qquad
\normalsize
\hbox{\rm then set}\quad
\small
B = \left[\begin{array}{cc}A^T A & C^T\\C&0 \end{array}\right]
\qquad r = \left[\begin{array}{c}0\\ \vdots\\1\end{array}\right] 
\normalsize
\end{equation}
to obtain the weights $\{a_k\}$ as the solution of $B a = r$.
For larger values of $m$, methods based on QR factorisation are advisable
but the cost of any aspect of the acceleration step is a very small
fraction of the cost of a direct solve of the generalised Westervelt equation.
The weights will typically be of mixed sign but an equally simple algorithm
allows the constraint $a_j \geq 0$ to be imposed.
For the examples involving the fixed point schemes 
\eqref{iteration_titr}, \eqref{iteration_fiti} this constraint gave weights
that were very close to the original Picard scheme, that is without any
effective acceleration.
As we will see shortly the unconstrained version gave 
considerable improvement to the fixed point schemes 
\eqref{iteration_titr}, \eqref{iteration_fiti}.

It has been shown that Anderson acceleration improves the convergence
rate of contractive fixed-point iterations in the vicinity of a fixed-point,
\cite{EvansPollockRebholzXiao:2020}, but will actually slow the rate of
quadratically convergent schemes such as those based on Newton methods.
The use of this acceleration method has so far been quite rare in the inverse problems literature, but see \cite{Jin_2020}.

\bigskip

Based on the analysis in Section \oldref{sec:forward}, we now investigate well-definedness of these iterative schemes.

\subsection{Well-definedness of projection based reconstruction schemes}

In the time trace data iteration scheme 
\begin{equation}\label{iteration_titr_k}
f_{k+1}(h(t))= 1-\frac{p_t(x_0,t;f_k))}{h'(t)}(1-f_k(p(x_0,t,f_k)))
\end{equation}
we would need $p_t\in L^2(0,T;C(\overline{\Omega}))$ to obtain a well-defined $L^2(0,T)$ trace $p_t(x_0,\cdot)$.
The final time data iteration scheme 
\begin{equation}\label{iteration_fiti_k}
f_{k+1}(g(x))= f_k(p(x,T;f_k)) + \frac{1}{p_t(x,T;f_k)} \, b\mathcal{A} (g(x)- p(x,T;f_k))\,,
\end{equation}
also requires evaluation of space derivatives.
To this end, we point out that by \cite[Lemma 3.3]{temam:2012}
we have $p=p(\cdot,\cdot;f)=G(f)\in C_w([0,T];H^2(\Omega))$.
Therefore, evaluation of $\mathcal{A} p(T)\in L^2(\Omega)$ makes sense.
However, since we also need $p_t$ to be at least in $C_w([0,T];L^\infty(\Omega))$ in order for the term 
$\frac{1}{p_t(\cdot,T)}$ to make sense and be bounded, we would need a higher order estimate here, which, however, would require higher differentiability of $f$.  

Since we prefer to stay with only Lipschitz continuous $f$ 
(bearing an mind that the inverse problem is less ill-posed on weaker regularity spaces)
we enforce well-definedness of the iteration schemes by projecting 
\begin{equation}\label{iteration_titr_proj}
f_{k+1}(h(t))= \mathcal{S}\left[1-\frac{\mathcal{P}_{[\underline{P},\overline{P}]}[p_t(x_0,t;f_k))]}{h'(t)}(1-f_k(p(x_0,t,f_k)))\right]
\end{equation}
\begin{equation}\label{iteration_fiti_proj}
f_{k+1}(g(x))= \mathcal{S}\left[f_k(p(x,T;f_k))  \mathcal{P}_{[\underline{Q},\overline{Q}]}\left[\frac{1}{p_t(x,T;f_k)}\right] \, b\mathcal{A} (g(x)- p(x,T;f_k))\right]\,.
\end{equation}
Here $\mathcal{P}_{[a,b]}$ is pointwise defined by $\mathcal{P}_{[a,b]}(z)=\max\{a,\max\{b,z\}\}$ and $\mathcal{S}$ is a smoothing operator mapping into $\mathcal{D}$ as defined in Corollary \oldref{cor:Fwelldef} and leaving functions already contained in $\mathcal{D}$ invariant, for example defining, for $y\in L^2(\Sigma)$, its smoothed version $\mathcal{S}[y]$ as a minimiser of
\[
\mbox{min} \|z-y\|_{L^2(\Sigma)}^2 \mbox{ such that } 
-\underline{\sigma}\leq f \leq \overline{\sigma} \,, \ -L\leq f' \leq L \mbox{ a.e. }
\] 

Here $\underline{P},\overline{P}$, $\underline{Q},\overline{Q}$, $\underline{\sigma},\overline{\sigma}$, $L$ should be chosen such that the exact solution $f_{\rm act}\in \mathcal{D}$ and 
$\underline{P}\leq p_t(x_0,t;f_{\rm act})\leq\overline{P}$ a.e. or 
$\underline{Q}\leq \frac{1}{p_t(x,T;f_{\rm act})}\leq\overline{Q}\;$ a.e., respectively.

\begin{corollary}\label{cor:itwelldef}
Under the conditions of Corollary \oldref{cor:Fwelldef}, the iteration schemes 
\begin{itemize}
\item \eqref{iteration_titr_proj} with $h'\in L^\infty(\Omega)$, $|h'|\geq\underline{\gamma}>0\;$ a.e., $\underline{P}\leq\overline{P}$
and 
\item \eqref{iteration_fiti_proj} with $g\in H^2(\omega)$, $|(g\circ \xi)'|\geq\underline{\gamma}>0$ a.e.,  for some curve $\xi$ contained in $\omega$, $\underline{Q}\leq\overline{Q}$ 
\end{itemize}
are well-defined.
\end{corollary}

\begin{remark}\label{rem:convfp}
Assume that the inverse problem has a solution $(f_{act},p_{act})$ such that $f_{act}\in\mathcal{D}$. 
Since $\mathcal{D}$ is bounded in $W^{1,\infty}(-M,M)$ we can conclude existence of a subsequnce $f_{k_n}$ of the projected fixed point iterates according to \eqref{iteration_titr_proj} or \eqref{iteration_fiti_proj}, respectively that converges weakly* to some $f^*\in\mathcal{D}$, which due to Proposition \oldref{prop:Fcont} (and its proof) together with $p^*=G(f)$ solves the inverse problem. 
In case the solution to the inverse problem is unique, a subsequence-subsequence argument yields weak* convergence in $W^{1,\infty}(-M,M)$ (hence norm convergence in $C([-M,M])$ and in $W^{s,\rho}(-M,M)$ for any $s<1$, $\rho\in[1,\infty)$) of the entire sequence of iterates.
\end{remark}

\subsection{Well-definedness of Newton's method}
The linearization of the forward operator $F=\mbox{tr}_\Sigma\circ G$ at some $f$ is defined by $F'(f)=\mbox{tr}_\Sigma\circ G'(f)$, where $G'(f)\underline{df}=\underline{dp}$ solves \eqref{Westervelt_f_PDEinit_lin}
that is, $\underline{dp}=\ppsi_t$ where $\ppsi$ solves \eqref{Westervelt_u_lin} with 
$\sigma=f(p)$, $\eta=f'(p)p_t$, $\tilde{r}=\underline{df} p_t$, $\ppsi_0=0$, $\ppsi_1=0$.
In order to obtain a well-defined $L^{\rm p}$ trace on $\Sigma$, we need at least 
$\underline{dp}=\ppsi_t\in L^{\rm p}([0,T];C(\Omega))$ in the time trace data and
$\underline{dp}=\ppsi_t\in C_w([0,T];L^{\rm p}(\Omega))$ in the final time data case, respectively.
We employ Lemma \oldref{lem:enest_lin} and the estimates
\[
\begin{aligned}
\|\eta\|_{L^2(0,T;L^q(\Omega))} 
&\leq L \|p_t\|_{L^2(0,T;L^q(\Omega))}\leq 
\begin{cases}
L C_{PF}^\Omega m_2 \mbox{ with }q=2 \ \mbox{ if }d=1\\
L C_{H^1,L^3}^\Omega m_2 \mbox{ with }q=3 \ \mbox{ if }d\in\{2,3\}
\end{cases} \\
\|\nabla\eta\|_{L^2(0,T;L^2(\Omega))} 
&\leq\|f''\|_{L^\infty(-M,M)} \|p_t\nabla p\|_{L^2(0,T;L^2(\Omega))}
+L\|\nabla p_t\|_{L^2(0,T;L^2(\Omega))}\\
&\leq\|f''\|_{L^\infty(-M,M)} (C_{H^1,L^4}^\Omega)^2 m_1 m_2
+L m_2
\end{aligned}
\]
for $f\in \mathcal{D}$, and $m_1$, $m_2$ as in \eqref{defM_mult_Autt}, where in the latter case we additionally have to assume $f\in W^{2,\infty}(-M,M)$. 
From this we deduce that the Newton iteration is well-defined by a minimiser  
\begin{equation}\label{eqn:Newton}
\begin{aligned}
f_{k+1}\in\ &\mbox{argmin}_{f_+\in\mathcal{D}} \|\mbox{tr}_\Sigma(p+\underline{dp})-y\|_{L^2(\Sigma)}^2  \\
&\mbox{such that $\;p$ solves \eqref{Westervelt_f_PDEinit} with $f=f_k$, and $\underline{dp}$ solves \eqref{Westervelt_f_PDEinit_lin} with $f=f_k$, $\underline{df}=f_+-f_k$}
\end{aligned}
\end{equation}
where $y=h$ in the time trace and $y=g$ in the final time data case.

\begin{corollary}\label{cor:Newtonwelldef}
Under the conditions of Corollary \oldref{cor:Fwelldef}, Newton's method is well-defined by \eqref{eqn:Newton} in the final time data case with $d\in\{1,2,3\}$ and in the time trace data case with $d=1$. The latter extends to $d\in\{2,3\}$ if $\mathcal{D}$ is replaced by $\tilde{\mathcal{D}}:=\mathcal{D}\cap\{f\in W^{2,\infty}(-M,M)\, : \, \|f''\|_{L^\infty(-M,M)}\leq N\}$ for some $N>0$.
\end{corollary}

A frozen version of \eqref{eqn:Newton} with fixed $f^0\in \mathcal{D}\cap W^{2,\infty}(-M,M)$, $p^0=G(f^0)$ can be defined by 
\begin{equation}\label{eqn:frozenNewton}
\begin{aligned}
f_{k+1}\in\ &\mbox{argmin}_{f_+\in\mathcal{D}} \|\mbox{tr}_\Sigma(p+\underline{dp})-y\|_{L^2(\Sigma)}^2  \\
&\mbox{such that $p$ solves \eqref{Westervelt_f_PDEinit} with $f=f_k$, and $\underline{dp}$ solves
 \eqref{Westervelt_f_PDEinit_lin} with $f=f^0$, $p=p^0$, $\underline{df}=f_+-f_k$.}
\end{aligned}
\end{equation}
This simplifies both the numerical computations and the analysis.
 
\begin{corollary}\label{cor:frozenNewtonwelldef}
For $f^0\in \mathcal{D}\cap W^{2,\infty}(-M,M)$, $p^0=G(f^0)$ a frozen Newton method is well-defined by \eqref{eqn:frozenNewton} in both data cases and space dimensions $d\in\{1,2,3\}$.
\end{corollary}

\begin{remark}\label{rem:convNewton}
Analogously to Remark \oldref{rem:convfp}, we can conclude subsequential weak* convergence of the  (frozen) Newton sequence to a solution of the inverse problem.
\end{remark}

\subsection*{Noisy data and regularisation}
In realistic measurements scenarios, the measured data $h$ or $g$ will typically be contaminated by noise and therefore only approximations $\tilde{h}\approx h$, $\tilde{g}\approx g$ are available. Their distance from the exact data can (if at all) only be estimated in some $L^{\rm p}$ norm corresponding to our choice of the data space $Y$, $\|\tilde{h}- h\|_Y\leq\delta$, $\|\tilde{g}- g\|_Y\leq\delta$, and they will typically lack differentiability. Still, the iteration schemes \eqref{iteration_titr_proj}, \eqref{iteration_fiti_proj} are well-defined upon replacement of $\tilde{h}$, $\tilde{g}$ by smoothed versions $\hat{h}$, $\hat{g}$ in such a way that the prerequisites of Corollary \oldref{cor:itwelldef} are preserved, that is, 
$\hat{h}'\in L^\infty(\Omega)$, $|\hat{h}'|\geq\underline{\hat{\gamma}}>0\;$ a.e., in the time trace case and $\hat{g}\in H^2(\omega)$, $|(\hat{g}\circ \xi)'|\geq\underline{\hat{\gamma}}>0$ a.e., for some curve $\xi$ contained in $\omega$, in the final time case. Indeed, if the exact data satisfies these conditions, the smoothed approximation can be chosen so that its distance from the exact data in the $L^{\rm p}$ norm is of the order of magnitude of the original noise level. This can for example be achieved by defining $\hat{d}=\hat{h}$ or $\hat{d}=\hat{g}\circ\xi$ from $\tilde{d}=\tilde{h}$ or $\tilde{d}=\tilde{g}\circ\xi$ according to 
\[
\hat{d}\in\mbox{argmin}_{b\in B_n:=\mbox{span}\{b_1,\ldots,b_n\}}  \|b-\tilde{d}\|_{L^{\rm p}} \mbox{ such that } s\, b' \geq \underline{\gamma}/2
\]
where $s = \mbox{sign}(d')$ and $\{b_1,\ldots,b_n\}$ is a sufficiently large set of smooth (in the final time case at least $H^2$) basis functions that allow to approximate the exact data up to precision $\mbox{dist}^{L^\infty}(d,B_n)\leq \underline{\gamma}/2$. Minimality of $\hat{d}$ and admissibility of $d$ for the above optimization problem immediately yield $\|\hat{d}-\tilde{d}\|_{L^{\rm p}}\leq \|d-\tilde{d}\|_{L^{\rm p}}$.
The optimization problem itself is well-posed and relatively simple to solve as a minimization of a convex cost function over a finite dimensional set under linear inequality constraints.
Clearly, this involves a priori information on the exact data, more precisely on the lower bound $\underline{\gamma}>0$ and the sign of $d'$, and additionally on the curve $\xi$ in the final time case. 

The regularization strategy that we employ here and in the fixed point schemes \eqref{iteration_titr_proj}, \eqref{iteration_fiti_proj} themselves as well as in the Newton type iterations \eqref{eqn:Newton}, \eqref{eqn:frozenNewton} is projection or restriction onto a (weakly) compact set. 
Indeed, considering a family of noisy data with $\delta\to0$ and stopping indices $k(\delta)\to\infty$, analogously to Remarks \oldref{rem:convfp}, \oldref{rem:convNewton} we obtain weak subsequential convergence of the iterates $f_{k(\delta)}$ to a solution of the inverse problem with exact data. Thus the proposed schemes are regularization methods in the classical sense of, e.g., \cite[Section 3]{EnglHankeNeubauer:1996}.
\section{Reconstructions}

In this section we show reconstructions of $f$ from either time trace or final time data.
The spatial set will be the interval $[0,1]$ and we will take the time trace measurement
point to be the right-hand endpoint $x=1$.

Our numerical implementation uses \eqref{Westervelt_f_PDEinit}
in the integrated version as in \eqref{Westervelt_k_int}
and so treat it as a parabolic equation with nonlocal memory
term $c^2\int_0^t \triangle p(\tau)\,d\tau$.
A Crank-Nicolson integrator was used with an inner iteration loop
to handle the nonlinear term $\,f(p)\,p_t$.
A Neumann boundary condition was imposed at the right hand endpoint;
the left hand condition could be Dirichlet, Neumann or impedance type.
Typically, in the physical model one would have zero initial conditions but
this isn't necessary for the mathematical formulation.

Data consisted of the measurements $h(t) = p(1,t)$ or $g(x) = p(x,T)$.
As a practical matter we used the above mentioned solver to obtain this data
and collected a sample at $50$ equally spaced points on the interval $[0,T]$ or $[0,1]$.
Uniformly distributed random noise was then added to these values to obtain
$h_{\rm meas}(t)$ or $g_{\rm meas}(x)$.
This was then pre-filtered by smoothing and up-resolving to the working
resolution of the number of points taken ($\sim$400 for the interval
$t\in [0,T]$ and $\sim$200 for the interval
$x\in [0,1]$) for the direct solver used in the inversion routine.

For the Newton scheme, the unknown $f$ was represented in terms of given basis functions.
Since we wish to make no constraints on the form of $f$,
we do not choose a basis with in-built restrictions as would be
obtained from an eigenfunction expansion.
Instead we used a sine basis that helped to realize the condition $f(0)=0$. 
In all cases the starting approximation for the iterative methods used
was the constant function $f\equiv0$.

We are going to show the results of the methods described in section \oldref{sec:reconstruction_schemes}. In particular, we will provide comparisons between frozen Newton and the fixed point iteration schemes from sections and \oldref{subsec:Newt} and  \oldref{subsec:Pic}, respectively. In the latter case, we will also show the effect of Anderson acceleration cf. section \oldref{subsec:And}. As to the parameters in the Anderson scheme, we simply fixed the depth to be $m=3$ and the weights $\beta_k=1$, which corresponds to full acceleration. 

\newbox\figurelegend
\newbox\figurelegendtwo
\newbox\figureone
\newbox\figuretwo
\newbox\figurethree
\newbox\figurefour
\newbox\figurefive
\newbox\figuresix
\newbox\figureseven
\newbox\figureeight
\newbox\figurenine
\newbox\figureten

\xfiglen=2.2 true in
\yfiglen=16.5 true in
\setbox\figureone=\vbox{\hsize=\xfiglen
\beginpicture
\eightrm
\footnotesize
  \setcoordinatesystem units <\xfiglen,\yfiglen> 
  \setplotarea x from 0 to 1, y from 0 to 0.08
  \axis bottom shiftedto y=0 ticks short numbered from 0 to 1 by 0.2 /
  \axis left ticks short numbered from 0 to 0.08 by 0.04 /
\footnotesize
\put {$u$} [rb] at 1 0.001
\put {$f(u)$} [l] at 0.01 0.08
\small
\linethickness=0.6pt
\setplotsymbol ({\sevenrm .})
\setquadratic
\setlinear
\setdashes
\Black{\relax  
\plot   
 0.0000   0.0000
  0.0200   0.0058
  0.0400   0.0113
  0.0600   0.0165
  0.0800   0.0213
  0.1000   0.0259
  0.1200   0.0301
  0.1400   0.0341
  0.1600   0.0379
  0.1800   0.0414
  0.2000   0.0447
  0.2200   0.0477
  0.2400   0.0506
  0.2600   0.0533
  0.2800   0.0557
  0.3000   0.0580
  0.3200   0.0602
  0.3400   0.0621
  0.3600   0.0640
  0.3800   0.0656
  0.4000   0.0672
  0.4200   0.0686
  0.4400   0.0698
  0.4600   0.0710
  0.4800   0.0720
  0.5000   0.0730
  0.5200   0.0738
  0.5400   0.0745
  0.5600   0.0752
  0.5800   0.0757
  0.6000   0.0762
  0.6200   0.0766
  0.6400   0.0769
  0.6600   0.0772
  0.6800   0.0773
  0.7000   0.0774
  0.7200   0.0775
  0.7400   0.0774
  0.7600   0.0774
  0.7800   0.0772
  0.8000   0.0770
  0.8200   0.0768
  0.8400   0.0765
  0.8600   0.0762
  0.8800   0.0758
  0.9000   0.0754
  0.9200   0.0749
  0.9400   0.0744
  0.9600   0.0739
  0.9800   0.0733
  1.0000   0.0727
/\relax}\relax
\setsolid
\Red{\relax
\plot
 0.0000   0.0000
  0.0200   0.0033
  0.0400   0.0061
  0.0600   0.0085
  0.0800   0.0105
  0.1000   0.0122
  0.1200   0.0137
  0.1400   0.0150
  0.1600   0.0161
  0.1800   0.0170
  0.2000   0.0177
  0.2200   0.0184
  0.2400   0.0189
  0.2600   0.0193
  0.2800   0.0196
  0.3000   0.0198
  0.3200   0.0200
  0.3400   0.0201
  0.3600   0.0201
  0.3800   0.0201
  0.4000   0.0200
  0.4200   0.0199
  0.4400   0.0198
  0.4600   0.0196
  0.4800   0.0195
  0.5000   0.0193
  0.5200   0.0190
  0.5400   0.0188
  0.5600   0.0186
  0.5800   0.0183
  0.6000   0.0181
  0.6200   0.0178
  0.6400   0.0175
  0.6600   0.0173
  0.6800   0.0170
  0.7000   0.0168
  0.7200   0.0165
  0.7400   0.0162
  0.7600   0.0160
  0.7800   0.0157
  0.8000   0.0155
  0.8200   0.0152
  0.8400   0.0150
  0.8600   0.0147
  0.8800   0.0145
  0.9000   0.0143
  0.9200   0.0140
  0.9400   0.0138
  0.9600   0.0136
  0.9800   0.0134
  1.0000   0.0132
 /\relax}\relax
\Blue{\relax
\plot
 0.0000   0.0000
  0.0200   0.0051
  0.0400   0.0093
  0.0600   0.0130
  0.0800   0.0161
  0.1000   0.0189
  0.1200   0.0212
  0.1400   0.0233
  0.1600   0.0250
  0.1800   0.0265
  0.2000   0.0278
  0.2200   0.0290
  0.2400   0.0299
  0.2600   0.0307
  0.2800   0.0313
  0.3000   0.0319
  0.3200   0.0323
  0.3400   0.0326
  0.3600   0.0329
  0.3800   0.0330
  0.4000   0.0331
  0.4200   0.0331
  0.4400   0.0331
  0.4600   0.0330
  0.4800   0.0329
  0.5000   0.0327
  0.5200   0.0325
  0.5400   0.0322
  0.5600   0.0320
  0.5800   0.0317
  0.6000   0.0314
  0.6200   0.0310
  0.6400   0.0307
  0.6600   0.0303
  0.6800   0.0299
  0.7000   0.0295
  0.7200   0.0292
  0.7400   0.0288
  0.7600   0.0284
  0.7800   0.0280
  0.8000   0.0276
  0.8200   0.0272
  0.8400   0.0268
  0.8600   0.0264
  0.8800   0.0260
  0.9000   0.0257
  0.9200   0.0253
  0.9400   0.0250
  0.9600   0.0246
  0.9800   0.0243
  1.0000   0.0240
/\relax}\relax
\Green{\relax
\plot
 -0.0000   0.0000
  0.0200   0.0060
  0.0400   0.0110
  0.0600   0.0154
  0.0800   0.0192
  0.1000   0.0225
  0.1200   0.0253
  0.1400   0.0278
  0.1600   0.0300
  0.1800   0.0320
  0.2000   0.0337
  0.2200   0.0351
  0.2400   0.0364
  0.2600   0.0376
  0.2800   0.0385
  0.3000   0.0393
  0.3200   0.0400
  0.3400   0.0406
  0.3600   0.0411
  0.3800   0.0415
  0.4000   0.0418
  0.4200   0.0420
  0.4400   0.0421
  0.4600   0.0422
  0.4800   0.0422
  0.5000   0.0421
  0.5200   0.0420
  0.5400   0.0418
  0.5600   0.0416
  0.5800   0.0414
  0.6000   0.0411
  0.6200   0.0408
  0.6400   0.0404
  0.6600   0.0401
  0.6800   0.0397
  0.7000   0.0393
  0.7200   0.0389
  0.7400   0.0384
  0.7600   0.0380
  0.7800   0.0375
  0.8000   0.0371
  0.8200   0.0366
  0.8400   0.0362
  0.8600   0.0357
  0.8800   0.0353
  0.9000   0.0348
  0.9200   0.0344
  0.9400   0.0340
  0.9600   0.0335
  0.9800   0.0331
  1.0000   0.0327
/\relax}\relax
\Orange{\relax
\plot
 0.0000   0.0000
  0.0200   0.0065
  0.0400   0.0120
  0.0600   0.0167
  0.0800   0.0208
  0.1000   0.0244
  0.1200   0.0276
  0.1400   0.0304
  0.1600   0.0329
  0.1800   0.0352
  0.2000   0.0371
  0.2200   0.0389
  0.2400   0.0405
  0.2600   0.0418
  0.2800   0.0431
  0.3000   0.0441
  0.3200   0.0451
  0.3400   0.0459
  0.3600   0.0466
  0.3800   0.0472
  0.4000   0.0477
  0.4200   0.0481
  0.4400   0.0484
  0.4600   0.0486
  0.4800   0.0488
  0.5000   0.0488
  0.5200   0.0488
  0.5400   0.0488
  0.5600   0.0487
  0.5800   0.0485
  0.6000   0.0483
  0.6200   0.0481
  0.6400   0.0478
  0.6600   0.0475
  0.6800   0.0471
  0.7000   0.0468
  0.7200   0.0464
  0.7400   0.0460
  0.7600   0.0455
  0.7800   0.0451
  0.8000   0.0446
  0.8200   0.0441
  0.8400   0.0437
  0.8600   0.0432
  0.8800   0.0427
  0.9000   0.0422
  0.9200   0.0417
  0.9400   0.0413
  0.9600   0.0408
  0.9800   0.0403
  1.0000   0.0398
/\relax}\relax
\Yellow{\relax
\plot
 0.0000   0.0000
  0.0200   0.0067
  0.0400   0.0125
  0.0600   0.0174
  0.0800   0.0218
  0.1000   0.0256
  0.1200   0.0290
  0.1400   0.0320
  0.1600   0.0347
  0.1800   0.0371
  0.2000   0.0393
  0.2200   0.0413
  0.2400   0.0431
  0.2600   0.0447
  0.2800   0.0461
  0.3000   0.0474
  0.3200   0.0485
  0.3400   0.0495
  0.3600   0.0504
  0.3800   0.0512
  0.4000   0.0519
  0.4200   0.0524
  0.4400   0.0529
  0.4600   0.0533
  0.4800   0.0536
  0.5000   0.0538
  0.5200   0.0539
  0.5400   0.0540
  0.5600   0.0540
  0.5800   0.0539
  0.6000   0.0538
  0.6200   0.0537
  0.6400   0.0535
  0.6600   0.0532
  0.6800   0.0529
  0.7000   0.0526
  0.7200   0.0522
  0.7400   0.0519
  0.7600   0.0515
  0.7800   0.0510
  0.8000   0.0506
  0.8200   0.0501
  0.8400   0.0497
  0.8600   0.0492
  0.8800   0.0487
  0.9000   0.0482
  0.9200   0.0477
  0.9400   0.0472
  0.9600   0.0467
  0.9800   0.0462
  1.0000   0.0457
/\relax}\relax
\Cyan{\relax
\plot
 0.0000   0.0000
  0.0200   0.0071
  0.0400   0.0131
  0.0600   0.0185
  0.0800   0.0232
  0.1000   0.0274
  0.1200   0.0312
  0.1400   0.0347
  0.1600   0.0379
  0.1800   0.0409
  0.2000   0.0436
  0.2200   0.0462
  0.2400   0.0485
  0.2600   0.0507
  0.2800   0.0528
  0.3000   0.0547
  0.3200   0.0564
  0.3400   0.0580
  0.3600   0.0594
  0.3800   0.0607
  0.4000   0.0619
  0.4200   0.0630
  0.4400   0.0640
  0.4600   0.0648
  0.4800   0.0656
  0.5000   0.0662
  0.5200   0.0667
  0.5400   0.0672
  0.5600   0.0676
  0.5800   0.0678
  0.6000   0.0681
  0.6200   0.0682
  0.6400   0.0683
  0.6600   0.0683
  0.6800   0.0683
  0.7000   0.0682
  0.7200   0.0681
  0.7400   0.0679
  0.7600   0.0677
  0.7800   0.0674
  0.8000   0.0671
  0.8200   0.0668
  0.8400   0.0664
  0.8600   0.0660
  0.8800   0.0656
  0.9000   0.0651
  0.9200   0.0646
  0.9400   0.0641
  0.9600   0.0636
  0.9800   0.0630
  1.0000   0.0624
/\relax}\relax
\Maroon{\relax
\plot
 0.0000   0.0000
  0.0200   0.0070
  0.0400   0.0130
  0.0600   0.0183
  0.0800   0.0231
  0.1000   0.0274
  0.1200   0.0314
  0.1400   0.0352
  0.1600   0.0387
  0.1800   0.0419
  0.2000   0.0450
  0.2200   0.0479
  0.2400   0.0507
  0.2600   0.0532
  0.2800   0.0556
  0.3000   0.0578
  0.3200   0.0599
  0.3400   0.0618
  0.3600   0.0635
  0.3800   0.0651
  0.4000   0.0666
  0.4200   0.0679
  0.4400   0.0691
  0.4600   0.0702
  0.4800   0.0712
  0.5000   0.0721
  0.5200   0.0728
  0.5400   0.0735
  0.5600   0.0741
  0.5800   0.0746
  0.6000   0.0750
  0.6200   0.0753
  0.6400   0.0756
  0.6600   0.0758
  0.6800   0.0759
  0.7000   0.0760
  0.7200   0.0760
  0.7400   0.0759
  0.7600   0.0758
  0.7800   0.0757
  0.8000   0.0755
  0.8200   0.0752
  0.8400   0.0749
  0.8600   0.0746
  0.8800   0.0742
  0.9000   0.0738
  0.9200   0.0733
  0.9400   0.0728
  0.9600   0.0723
  0.9800   0.0717
  1.0000   0.0711
/\relax}\relax
\endpicture
}
\setbox\figuretwo=\vbox{\hsize=\xfiglen
\beginpicture
\eightrm
\footnotesize
  \setcoordinatesystem units <\xfiglen,\yfiglen> 
  \setplotarea x from 0 to 1, y from 0 to 0.08
  \axis bottom shiftedto y=0 ticks short numbered from 0 to 1 by 0.2 /
  \axis left ticks short numbered from 0 to 0.08 by 0.04 /
\footnotesize
\put {$u$} [rb] at 1 0.001
\put {$f(u)$} [l] at 0.01 0.08
\small
\linethickness=0.6pt
\setplotsymbol ({\sevenrm .})
\setquadratic
\setlinear
\setdashes
\Black{\relax  
\plot   
 0.0000   0.0000
  0.0200   0.0058
  0.0400   0.0113
  0.0600   0.0165
  0.0800   0.0213
  0.1000   0.0259
  0.1200   0.0301
  0.1400   0.0341
  0.1600   0.0379
  0.1800   0.0414
  0.2000   0.0447
  0.2200   0.0477
  0.2400   0.0506
  0.2600   0.0533
  0.2800   0.0557
  0.3000   0.0580
  0.3200   0.0602
  0.3400   0.0621
  0.3600   0.0640
  0.3800   0.0656
  0.4000   0.0672
  0.4200   0.0686
  0.4400   0.0698
  0.4600   0.0710
  0.4800   0.0720
  0.5000   0.0730
  0.5200   0.0738
  0.5400   0.0745
  0.5600   0.0752
  0.5800   0.0757
  0.6000   0.0762
  0.6200   0.0766
  0.6400   0.0769
  0.6600   0.0772
  0.6800   0.0773
  0.7000   0.0774
  0.7200   0.0775
  0.7400   0.0774
  0.7600   0.0774
  0.7800   0.0772
  0.8000   0.0770
  0.8200   0.0768
  0.8400   0.0765
  0.8600   0.0762
  0.8800   0.0758
  0.9000   0.0754
  0.9200   0.0749
  0.9400   0.0744
  0.9600   0.0739
  0.9800   0.0733
  1.0000   0.0727
/\relax}\relax
\setsolid
\Red{\relax
\plot
 -0.0000   0.0000
  0.0200   0.0091
  0.0400   0.0168
  0.0600   0.0234
  0.0800   0.0292
  0.1000   0.0341
  0.1200   0.0385
  0.1400   0.0423
  0.1600   0.0456
  0.1800   0.0486
  0.2000   0.0511
  0.2200   0.0534
  0.2400   0.0553
  0.2600   0.0570
  0.2800   0.0584
  0.3000   0.0597
  0.3200   0.0607
  0.3400   0.0616
  0.3600   0.0623
  0.3800   0.0628
  0.4000   0.0632
  0.4200   0.0635
  0.4400   0.0637
  0.4600   0.0638
  0.4800   0.0637
  0.5000   0.0636
  0.5200   0.0634
  0.5400   0.0631
  0.5600   0.0628
  0.5800   0.0624
  0.6000   0.0619
  0.6200   0.0614
  0.6400   0.0609
  0.6600   0.0603
  0.6800   0.0596
  0.7000   0.0590
  0.7200   0.0583
  0.7400   0.0576
  0.7600   0.0569
  0.7800   0.0562
  0.8000   0.0555
  0.8200   0.0547
  0.8400   0.0540
  0.8600   0.0533
  0.8800   0.0526
  0.9000   0.0519
  0.9200   0.0513
  0.9400   0.0506
  0.9600   0.0500
  0.9800   0.0493
  1.0000   0.0487
/\relax}\relax
\Blue{\relax
\plot
 0.0000   0.0000
  0.0200   0.0065
  0.0400   0.0121
  0.0600   0.0170
  0.0800   0.0214
  0.1000   0.0253
  0.1200   0.0289
  0.1400   0.0322
  0.1600   0.0353
  0.1800   0.0382
  0.2000   0.0409
  0.2200   0.0435
  0.2400   0.0459
  0.2600   0.0482
  0.2800   0.0504
  0.3000   0.0524
  0.3200   0.0543
  0.3400   0.0561
  0.3600   0.0578
  0.3800   0.0593
  0.4000   0.0607
  0.4200   0.0620
  0.4400   0.0632
  0.4600   0.0642
  0.4800   0.0652
  0.5000   0.0660
  0.5200   0.0667
  0.5400   0.0672
  0.5600   0.0677
  0.5800   0.0681
  0.6000   0.0684
  0.6200   0.0685
  0.6400   0.0686
  0.6600   0.0687
  0.6800   0.0686
  0.7000   0.0685
  0.7200   0.0683
  0.7400   0.0681
  0.7600   0.0678
  0.7800   0.0674
  0.8000   0.0671
  0.8200   0.0667
  0.8400   0.0662
  0.8600   0.0657
  0.8800   0.0652
  0.9000   0.0647
  0.9200   0.0641
  0.9400   0.0636
  0.9600   0.0630
  0.9800   0.0623
  1.0000   0.0617
/\relax}\relax
\Green{\relax
\plot
 0.0000   0.0000
  0.0200   0.0073
  0.0400   0.0135
  0.0600   0.0190
  0.0800   0.0239
  0.1000   0.0283
  0.1200   0.0324
  0.1400   0.0361
  0.1600   0.0396
  0.1800   0.0429
  0.2000   0.0459
  0.2200   0.0488
  0.2400   0.0515
  0.2600   0.0540
  0.2800   0.0564
  0.3000   0.0586
  0.3200   0.0606
  0.3400   0.0626
  0.3600   0.0643
  0.3800   0.0660
  0.4000   0.0675
  0.4200   0.0688
  0.4400   0.0701
  0.4600   0.0712
  0.4800   0.0722
  0.5000   0.0730
  0.5200   0.0738
  0.5400   0.0745
  0.5600   0.0751
  0.5800   0.0756
  0.6000   0.0760
  0.6200   0.0763
  0.6400   0.0766
  0.6600   0.0768
  0.6800   0.0769
  0.7000   0.0770
  0.7200   0.0770
  0.7400   0.0770
  0.7600   0.0769
  0.7800   0.0767
  0.8000   0.0765
  0.8200   0.0762
  0.8400   0.0759
  0.8600   0.0755
  0.8800   0.0751
  0.9000   0.0747
  0.9200   0.0742
  0.9400   0.0736
  0.9600   0.0730
  0.9800   0.0724
  1.0000   0.0718
/\relax}\relax
\endpicture
}
\setbox\figurethree=\vbox{\hsize=\xfiglen
\beginpicture
\eightrm
\footnotesize
  \setcoordinatesystem units <\xfiglen,\yfiglen> 
  \setplotarea x from 0 to 1, y from 0 to 0.08
  \axis bottom shiftedto y=0 ticks short numbered from 0 to 1 by 0.2 /
  \axis left ticks short numbered from 0 to 0.08 by 0.04 /
\footnotesize
\put {$u$} [rb] at 1 0.001
\put {$f(u)$} [l] at 0.01 0.08
\small
\linethickness=0.6pt
\setplotsymbol ({\sevenrm .})
\setquadratic
\setlinear
\setdashes
\Black{\relax  
\plot   
  0.0000   0.0000
  0.0200   0.0058
  0.0400   0.0113
  0.0600   0.0165
  0.0800   0.0213
  0.1000   0.0259
  0.1200   0.0301
  0.1400   0.0341
  0.1600   0.0379
  0.1800   0.0414
  0.2000   0.0447
  0.2200   0.0477
  0.2400   0.0506
  0.2600   0.0533
  0.2800   0.0557
  0.3000   0.0580
  0.3200   0.0602
  0.3400   0.0621
  0.3600   0.0640
  0.3800   0.0656
  0.4000   0.0672
  0.4200   0.0686
  0.4400   0.0698
  0.4600   0.0710
  0.4800   0.0720
  0.5000   0.0730
  0.5200   0.0738
  0.5400   0.0745
  0.5600   0.0752
  0.5800   0.0757
  0.6000   0.0762
  0.6200   0.0766
  0.6400   0.0769
  0.6600   0.0772
  0.6800   0.0773
  0.7000   0.0774
  0.7200   0.0775
  0.7400   0.0774
  0.7600   0.0774
  0.7800   0.0772
  0.8000   0.0770
  0.8200   0.0768
  0.8400   0.0765
  0.8600   0.0762
  0.8800   0.0758
  0.9000   0.0754
  0.9200   0.0749
  0.9400   0.0744
  0.9600   0.0739
  0.9800   0.0733
  1.0000   0.0727
/\relax}\relax
\setsolid
\Red{\relax
\plot
 0.0000  -0.0000
  0.0200   0.0057
  0.0400   0.0113
  0.0600   0.0167
  0.0800   0.0220
  0.1000   0.0271
  0.1200   0.0318
  0.1400   0.0363
  0.1600   0.0405
  0.1800   0.0443
  0.2000   0.0479
  0.2200   0.0512
  0.2400   0.0542
  0.2600   0.0571
  0.2800   0.0597
  0.3000   0.0621
  0.3200   0.0644
  0.3400   0.0665
  0.3600   0.0685
  0.3800   0.0703
  0.4000   0.0720
  0.4200   0.0735
  0.4400   0.0749
  0.4600   0.0761
  0.4800   0.0772
  0.5000   0.0782
  0.5200   0.0790
  0.5400   0.0797
  0.5600   0.0803
  0.5800   0.0808
  0.6000   0.0812
  0.6200   0.0815
  0.6400   0.0818
  0.6600   0.0820
  0.6800   0.0821
  0.7000   0.0821
  0.7200   0.0821
  0.7400   0.0820
  0.7600   0.0818
  0.7800   0.0815
  0.8000   0.0811
  0.8200   0.0807
  0.8400   0.0803
  0.8600   0.0798
  0.8800   0.0793
  0.9000   0.0787
  0.9200   0.0781
  0.9400   0.0774
  0.9600   0.0765
  0.9800   0.0755
  1.0000   0.0742
/\relax}\relax
\Blue{\relax
\plot
 0.0000  -0.0000
  0.0200   0.0055
  0.0400   0.0110
  0.0600   0.0163
  0.0800   0.0214
  0.1000   0.0262
  0.1200   0.0307
  0.1400   0.0348
  0.1600   0.0387
  0.1800   0.0422
  0.2000   0.0454
  0.2200   0.0483
  0.2400   0.0509
  0.2600   0.0534
  0.2800   0.0556
  0.3000   0.0577
  0.3200   0.0597
  0.3400   0.0616
  0.3600   0.0633
  0.3800   0.0650
  0.4000   0.0665
  0.4200   0.0680
  0.4400   0.0693
  0.4600   0.0705
  0.4800   0.0715
  0.5000   0.0725
  0.5200   0.0733
  0.5400   0.0740
  0.5600   0.0745
  0.5800   0.0750
  0.6000   0.0755
  0.6200   0.0758
  0.6400   0.0761
  0.6600   0.0764
  0.6800   0.0766
  0.7000   0.0768
  0.7200   0.0769
  0.7400   0.0769
  0.7600   0.0769
  0.7800   0.0768
  0.8000   0.0766
  0.8200   0.0764
  0.8400   0.0761
  0.8600   0.0758
  0.8800   0.0755
  0.9000   0.0752
  0.9200   0.0750
  0.9400   0.0748
  0.9600   0.0746
  0.9800   0.0744
  1.0000   0.0742
/\relax}\relax
\Green{\relax
\plot
  0.0000  -0.0000
  0.0200   0.0055
  0.0400   0.0109
  0.0600   0.0161
  0.0800   0.0212
  0.1000   0.0260
  0.1200   0.0306
  0.1400   0.0348
  0.1600   0.0387
  0.1800   0.0423
  0.2000   0.0456
  0.2200   0.0487
  0.2400   0.0514
  0.2600   0.0540
  0.2800   0.0563
  0.3000   0.0584
  0.3200   0.0604
  0.3400   0.0623
  0.3600   0.0640
  0.3800   0.0657
  0.4000   0.0672
  0.4200   0.0686
  0.4400   0.0699
  0.4600   0.0711
  0.4800   0.0722
  0.5000   0.0731
  0.5200   0.0740
  0.5400   0.0747
  0.5600   0.0753
  0.5800   0.0758
  0.6000   0.0762
  0.6200   0.0765
  0.6400   0.0768
  0.6600   0.0770
  0.6800   0.0771
  0.7000   0.0772
  0.7200   0.0773
  0.7400   0.0773
  0.7600   0.0772
  0.7800   0.0771
  0.8000   0.0770
  0.8200   0.0768
  0.8400   0.0765
  0.8600   0.0762
  0.8800   0.0759
  0.9000   0.0755
  0.9200   0.0751
  0.9400   0.0746
  0.9600   0.0741
  0.9800   0.0734
  1.0000   0.0726
/\relax}\relax
\endpicture
}

\xfiglen=2.4 true in
\yfiglen=3.8 true in
\setbox\figurefour=\vbox{\hsize=\xfiglen
\beginpicture
\eightrm
\footnotesize
  \setcoordinatesystem units <\xfiglen,\yfiglen>  point at  0 -0.20
  \setplotarea x from 0 to 1, y from -0.2 to 0.2
  \axis bottom shiftedto y=0 ticks short numbered from 0.2 to 1 by 0.2 /
  \axis left ticks short numbered from -0.1 to 0.2 by 0.1 /
\footnotesize
\put {$u$} [rb] at 1 0.01
\put {$f(u)$} [l] at 0.01 0.18
\small
\linethickness=0.6pt
\setplotsymbol ({\sevenrm .})
\setquadratic
\setlinear
\setdashes
\Black{\relax  
\plot   
  0.0000   0.0000
  0.0200   0.0083
  0.0400   0.0292
  0.0600   0.0547
  0.0800   0.0761
  0.1000   0.0858
  0.1200   0.0791
  0.1400   0.0548
  0.1600   0.0155
  0.1800  -0.0328
  0.2000  -0.0821
  0.2200  -0.1238
  0.2400  -0.1499
  0.2600  -0.1552
  0.2800  -0.1374
  0.3000  -0.0981
  0.3200  -0.0424
  0.3400   0.0218
  0.3600   0.0852
  0.3800   0.1383
  0.4000   0.1730
  0.4200   0.1835
  0.4400   0.1681
  0.4600   0.1282
  0.4800   0.0695
  0.5000   0.0000
  0.5200  -0.0704
  0.5400  -0.1315
  0.5600  -0.1747
  0.5800  -0.1935
  0.6000  -0.1850
  0.6200  -0.1504
  0.6400  -0.0943
  0.6600  -0.0246
  0.6800   0.0489
  0.7000   0.1158
  0.7200   0.1666
  0.7400   0.1941
  0.7600   0.1944
  0.7800   0.1673
  0.8000   0.1166
  0.8200   0.0494
  0.8400  -0.0249
  0.8600  -0.0958
  0.8800  -0.1533
  0.9000  -0.1894
  0.9200  -0.1988
  0.9400  -0.1803
  0.9600  -0.1365
  0.9800  -0.0734
  1.0000  -0.0000
/\relax}\relax
\setsolid
\Red{\relax 
\plot
  0.0000   0.0000
  0.0200   0.0083
  0.0400   0.0292
  0.0600   0.0545
  0.0800   0.0757
  0.1000   0.0850
  0.1200   0.0779
  0.1400   0.0533
  0.1600   0.0139
  0.1800  -0.0342
  0.2000  -0.0829
  0.2200  -0.1234
  0.2400  -0.1482
  0.2600  -0.1520
  0.2800  -0.1328
  0.3000  -0.0925
  0.3200  -0.0365
  0.3400   0.0273
  0.3600   0.0895
  0.3800   0.1405
  0.4000   0.1724
  0.4200   0.1802
  0.4400   0.1620
  0.4600   0.1203
  0.4800   0.0606
  0.5000  -0.0086
  0.5200  -0.0773
  0.5400  -0.1356
  0.5600  -0.1750
  0.5800  -0.1897
  0.6000  -0.1775
  0.6200  -0.1400
  0.6400  -0.0824
  0.6600  -0.0130
  0.6800   0.0584
  0.7000   0.1216
  0.7200   0.1676
  0.7400   0.1897
  0.7600   0.1847
  0.7800   0.1533
  0.8000   0.0999
  0.8200   0.0321
  0.8400  -0.0406
  0.8600  -0.1078
  0.8800  -0.1600
  0.9000  -0.1899
  0.9200  -0.1933
  0.9400  -0.1697
  0.9600  -0.1226
  0.9800  -0.0586
  1.0000   0.0132
/\relax}\relax
\Blue{\relax 
\plot
  0.0000   0.0000
  0.0200   0.0084
  0.0400   0.0292
  0.0600   0.0547
  0.0800   0.0761
  0.1000   0.0858
  0.1200   0.0790
  0.1400   0.0547
  0.1600   0.0154
  0.1800  -0.0329
  0.2000  -0.0822
  0.2200  -0.1238
  0.2400  -0.1499
  0.2600  -0.1551
  0.2800  -0.1373
  0.3000  -0.0980
  0.3200  -0.0422
  0.3400   0.0221
  0.3600   0.0856
  0.3800   0.1387
  0.4000   0.1733
  0.4200   0.1838
  0.4400   0.1682
  0.4600   0.1282
  0.4800   0.0693
  0.5000  -0.0003
  0.5200  -0.0707
  0.5400  -0.1319
  0.5600  -0.1750
  0.5800  -0.1937
  0.6000  -0.1851
  0.6200  -0.1503
  0.6400  -0.0939
  0.6600  -0.0240
  0.6800   0.0498
  0.7000   0.1168
  0.7200   0.1677
  0.7400   0.1950
  0.7600   0.1948
  0.7800   0.1672
  0.8000   0.1160
  0.8200   0.0484
  0.8400  -0.0262
  0.8600  -0.0972
  0.8800  -0.1547
  0.9000  -0.1907
  0.9200  -0.1999
  0.9400  -0.1811
  0.9600  -0.1368
  0.9800  -0.0733
  1.0000   0.0006
/\relax}\relax
\Green{\relax 
\plot
  0.0000   0.0000
  0.0200   0.0084
  0.0400   0.0292
  0.0600   0.0547
  0.0800   0.0761
  0.1000   0.0858
  0.1200   0.0791
  0.1400   0.0547
  0.1600   0.0155
  0.1800  -0.0328
  0.2000  -0.0821
  0.2200  -0.1237
  0.2400  -0.1499
  0.2600  -0.1551
  0.2800  -0.1374
  0.3000  -0.0981
  0.3200  -0.0424
  0.3400   0.0218
  0.3600   0.0852
  0.3800   0.1383
  0.4000   0.1729
  0.4200   0.1835
  0.4400   0.1680
  0.4600   0.1282
  0.4800   0.0695
  0.5000   0.0000
  0.5200  -0.0703
  0.5400  -0.1315
  0.5600  -0.1746
  0.5800  -0.1934
  0.6000  -0.1850
  0.6200  -0.1503
  0.6400  -0.0943
  0.6600  -0.0246
  0.6800   0.0489
  0.7000   0.1158
  0.7200   0.1667
  0.7400   0.1942
  0.7600   0.1945
  0.7800   0.1674
  0.8000   0.1166
  0.8200   0.0494
  0.8400  -0.0249
  0.8600  -0.0958
  0.8800  -0.1533
  0.9000  -0.1894
  0.9200  -0.1989
  0.9400  -0.1804
  0.9600  -0.1366
  0.9800  -0.0735
  1.0000  -0.0001
/\relax}\relax
%
\endpicture
}

\setbox\figurefive=\vbox{\hsize=\xfiglen
\beginpicture
\eightrm
\footnotesize
  \setcoordinatesystem units <\xfiglen,\yfiglen>  point at  0 -0.20
  \setplotarea x from 0 to 1, y from -0.2 to 0.2
  \axis bottom shiftedto y=0 ticks short numbered from 0.2 to 1 by 0.2 /
  \axis left ticks short numbered from 0.0 to 0.2 by 0.1 /
\footnotesize
\put {$u$} [rb] at 1 0.01
\put {$f(u)$} [l] at 0.01 0.18
\small
\linethickness=0.6pt
\setplotsymbol ({\sevenrm .})
\setquadratic
\setlinear
\setdashes
\Black{\relax  
\plot   
  0.0000   0.0000
  0.0200   0.0083
  0.0400   0.0292
  0.0600   0.0547
  0.0800   0.0761
  0.1000   0.0858
  0.1200   0.0791
  0.1400   0.0548
  0.1600   0.0155
  0.1800  -0.0328
  0.2000  -0.0821
  0.2200  -0.1238
  0.2400  -0.1499
  0.2600  -0.1552
  0.2800  -0.1374
  0.3000  -0.0981
  0.3200  -0.0424
  0.3400   0.0218
  0.3600   0.0852
  0.3800   0.1383
  0.4000   0.1730
  0.4200   0.1835
  0.4400   0.1681
  0.4600   0.1282
  0.4800   0.0695
  0.5000   0.0000
  0.5200  -0.0704
  0.5400  -0.1315
  0.5600  -0.1747
  0.5800  -0.1935
  0.6000  -0.1850
  0.6200  -0.1504
  0.6400  -0.0943
  0.6600  -0.0246
  0.6800   0.0489
  0.7000   0.1158
  0.7200   0.1666
  0.7400   0.1941
  0.7600   0.1944
  0.7800   0.1673
  0.8000   0.1166
  0.8200   0.0494
  0.8400  -0.0249
  0.8600  -0.0958
  0.8800  -0.1533
  0.9000  -0.1894
  0.9200  -0.1988
  0.9400  -0.1803
  0.9600  -0.1365
  0.9800  -0.0734
  1.0000  -0.0000
/\relax}\relax
\setsolid
\Red{\relax 
\plot
  0.0000   0.0000
  0.0200   0.0084
  0.0400   0.0292
  0.0600   0.0547
  0.0800   0.0761
  0.1000   0.0858
  0.1200   0.0791
  0.1400   0.0547
  0.1600   0.0155
  0.1800  -0.0329
  0.2000  -0.0822
  0.2200  -0.1238
  0.2400  -0.1499
  0.2600  -0.1552
  0.2800  -0.1374
  0.3000  -0.0981
  0.3200  -0.0424
  0.3400   0.0220
  0.3600   0.0855
  0.3800   0.1387
  0.4000   0.1733
  0.4200   0.1839
  0.4400   0.1683
  0.4600   0.1284
  0.4800   0.0695
  0.5000  -0.0001
  0.5200  -0.0706
  0.5400  -0.1319
  0.5600  -0.1750
  0.5800  -0.1938
  0.6000  -0.1853
  0.6200  -0.1505
  0.6400  -0.0942
  0.6600  -0.0242
  0.6800   0.0496
  0.7000   0.1167
  0.7200   0.1677
  0.7400   0.1951
  0.7600   0.1951
  0.7800   0.1675
  0.8000   0.1164
  0.8200   0.0487
  0.8400  -0.0259
  0.8600  -0.0970
  0.8800  -0.1546
  0.9000  -0.1907
  0.9200  -0.2000
  0.9400  -0.1813
  0.9600  -0.1371
  0.9800  -0.0736
  1.0000   0.0003
/\relax}\relax
\Blue{\relax 
\plot
  0.0000   0.0000
  0.0200   0.0084
  0.0400   0.0292
  0.0600   0.0547
  0.0800   0.0761
  0.1000   0.0858
  0.1200   0.0791
  0.1400   0.0547
  0.1600   0.0154
  0.1800  -0.0329
  0.2000  -0.0821
  0.2200  -0.1237
  0.2400  -0.1499
  0.2600  -0.1551
  0.2800  -0.1373
  0.3000  -0.0981
  0.3200  -0.0424
  0.3400   0.0218
  0.3600   0.0852
  0.3800   0.1383
  0.4000   0.1729
  0.4200   0.1835
  0.4400   0.1680
  0.4600   0.1282
  0.4800   0.0695
  0.5000  -0.0000
  0.5200  -0.0704
  0.5400  -0.1315
  0.5600  -0.1746
  0.5800  -0.1934
  0.6000  -0.1849
  0.6200  -0.1503
  0.6400  -0.0942
  0.6600  -0.0246
  0.6800   0.0489
  0.7000   0.1158
  0.7200   0.1666
  0.7400   0.1942
  0.7600   0.1944
  0.7800   0.1673
  0.8000   0.1166
  0.8200   0.0494
  0.8400  -0.0249
  0.8600  -0.0958
  0.8800  -0.1533
  0.9000  -0.1893
  0.9200  -0.1988
  0.9400  -0.1803
  0.9600  -0.1365
  0.9800  -0.0735
  1.0000  -0.0001
/\relax}\relax
%
%
\endpicture
}

\setbox\figuresix=\vbox{\hsize=\xfiglen
\beginpicture
\eightrm
\footnotesize
  \setcoordinatesystem units <\xfiglen,\yfiglen>  point at  0 -0.20
  \setplotarea x from 0 to 1, y from -0.2 to 0.2
  \axis bottom shiftedto y=0 ticks short numbered from 0.2 to 1 by 0.2 /
  \axis left ticks short numbered from 0.0 to 0.2 by 0.1 /
\footnotesize
\put {$u$} [rb] at 1 0.01
\put {$f(u)$} [l] at 0.01 0.18
\small
\linethickness=0.6pt
\setplotsymbol ({\sevenrm .})
\setquadratic
\setlinear
\setdashes
\Black{\relax  
\plot   
  0.0000   0.0000
  0.0200   0.0083
  0.0400   0.0292
  0.0600   0.0547
  0.0800   0.0761
  0.1000   0.0858
  0.1200   0.0791
  0.1400   0.0548
  0.1600   0.0155
  0.1800  -0.0328
  0.2000  -0.0821
  0.2200  -0.1238
  0.2400  -0.1499
  0.2600  -0.1552
  0.2800  -0.1374
  0.3000  -0.0981
  0.3200  -0.0424
  0.3400   0.0218
  0.3600   0.0852
  0.3800   0.1383
  0.4000   0.1730
  0.4200   0.1835
  0.4400   0.1681
  0.4600   0.1282
  0.4800   0.0695
  0.5000   0.0000
  0.5200  -0.0704
  0.5400  -0.1315
  0.5600  -0.1747
  0.5800  -0.1935
  0.6000  -0.1850
  0.6200  -0.1504
  0.6400  -0.0943
  0.6600  -0.0246
  0.6800   0.0489
  0.7000   0.1158
  0.7200   0.1666
  0.7400   0.1941
  0.7600   0.1944
  0.7800   0.1673
  0.8000   0.1166
  0.8200   0.0494
  0.8400  -0.0249
  0.8600  -0.0958
  0.8800  -0.1533
  0.9000  -0.1894
  0.9200  -0.1988
  0.9400  -0.1803
  0.9600  -0.1365
  0.9800  -0.0734
  1.0000  -0.0000
/\relax}\relax
\setsolid
\Red{\relax 
\plot
  0.0000   0.0000
  0.0200   0.0090
  0.0400   0.0292
  0.0600   0.0594
  0.0800   0.0859
  0.1000   0.0940
  0.1200   0.0780
  0.1400   0.0428
  0.1600  -0.0023
  0.1800  -0.0482
  0.2000  -0.0877
  0.2200  -0.1162
  0.2400  -0.1316
  0.2600  -0.1341
  0.2800  -0.1244
  0.3000  -0.1014
  0.3200  -0.0622
  0.3400  -0.0046
  0.3600   0.0688
  0.3800   0.1459
  0.4000   0.2046
  0.4200   0.2216
  0.4400   0.1873
  0.4600   0.1143
  0.4800   0.0311
  0.5000  -0.0387
  0.5200  -0.0890
  0.5400  -0.1260
  0.5600  -0.1533
  0.5800  -0.1658
  0.6000  -0.1591
  0.6200  -0.1363
  0.6400  -0.1030
  0.6600  -0.0557
  0.6800   0.0144
  0.7000   0.1049
  0.7200   0.1912
  0.7400   0.2393
  0.7600   0.2307
  0.7800   0.1740
  0.8000   0.0936
  0.8200   0.0117
  0.8400  -0.0583
  0.8600  -0.1094
  0.8800  -0.1393
  0.9000  -0.1537
  0.9200  -0.1614
  0.9400  -0.1622
  0.9600  -0.1415
  0.9800  -0.0852
  1.0000  -0.0000
/\relax}\relax
\Blue{\relax 
\plot
  0.0000   0.0000
  0.0200   0.0100
  0.0400   0.0272
  0.0600   0.0514
  0.0800   0.0749
  0.1000   0.0878
  0.1200   0.0833
  0.1400   0.0586
  0.1600   0.0158
  0.1800  -0.0372
  0.2000  -0.0880
  0.2200  -0.1251
  0.2400  -0.1447
  0.2600  -0.1491
  0.2800  -0.1382
  0.3000  -0.1061
  0.3200  -0.0489
  0.3400   0.0245
  0.3600   0.0934
  0.3800   0.1404
  0.4000   0.1647
  0.4200   0.1752
  0.4400   0.1731
  0.4600   0.1458
  0.4800   0.0824
  0.5000  -0.0064
  0.5200  -0.0905
  0.5400  -0.1446
  0.5600  -0.1678
  0.5800  -0.1762
  0.6000  -0.1779
  0.6200  -0.1608
  0.6400  -0.1093
  0.6600  -0.0273
  0.6800   0.0591
  0.7000   0.1238
  0.7200   0.1611
  0.7400   0.1817
  0.7600   0.1915
  0.7800   0.1804
  0.8000   0.1341
  0.8200   0.0546
  0.8400  -0.0360
  0.8600  -0.1111
  0.8800  -0.1580
  0.9000  -0.1808
  0.9200  -0.1873
  0.9400  -0.1773
  0.9600  -0.1433
  0.9800  -0.0813
  1.0000  -0.0000
/\relax}\relax
\Green{\relax 
\plot
  0.0000   0.0000
  0.0200   0.0098
  0.0400   0.0279
  0.0600   0.0533
  0.0800   0.0762
  0.1000   0.0861
  0.1200   0.0785
  0.1400   0.0545
  0.1600   0.0167
  0.1800  -0.0313
  0.2000  -0.0822
  0.2200  -0.1253
  0.2400  -0.1510
  0.2600  -0.1547
  0.2800  -0.1366
  0.3000  -0.0987
  0.3200  -0.0436
  0.3400   0.0223
  0.3600   0.0878
  0.3800   0.1403
  0.4000   0.1713
  0.4200   0.1789
  0.4400   0.1647
  0.4600   0.1295
  0.4800   0.0746
  0.5000   0.0045
  0.5200  -0.0699
  0.5400  -0.1347
  0.5600  -0.1781
  0.5800  -0.1945
  0.6000  -0.1841
  0.6200  -0.1496
  0.6400  -0.0948
  0.6600  -0.0250
  0.6800   0.0505
  0.7000   0.1187
  0.7200   0.1676
  0.7400   0.1910
  0.7600   0.1892
  0.7800   0.1649
  0.8000   0.1197
  0.8200   0.0551
  0.8400  -0.0220
  0.8600  -0.0978
  0.8800  -0.1571
  0.9000  -0.1904
  0.9200  -0.1969
  0.9400  -0.1790
  0.9600  -0.1380
  0.9800  -0.0758
  1.0000  -0.0000
/\relax}\relax
\Orange{\relax 
\plot
  0.0000   0.0000
  0.0200   0.0097
  0.0400   0.0276
  0.0600   0.0532
  0.0800   0.0767
  0.1000   0.0870
  0.1200   0.0789
  0.1400   0.0535
  0.1600   0.0150
  0.1800  -0.0320
  0.2000  -0.0812
  0.2200  -0.1238
  0.2400  -0.1507
  0.2600  -0.1555
  0.2800  -0.1370
  0.3000  -0.0976
  0.3200  -0.0425
  0.3400   0.0215
  0.3600   0.0852
  0.3800   0.1387
  0.4000   0.1734
  0.4200   0.1835
  0.4400   0.1672
  0.4600   0.1269
  0.4800   0.0688
  0.5000   0.0007
  0.5200  -0.0687
  0.5400  -0.1303
  0.5600  -0.1751
  0.5800  -0.1949
  0.6000  -0.1858
  0.6200  -0.1497
  0.6400  -0.0931
  0.6600  -0.0244
  0.6800   0.0480
  0.7000   0.1152
  0.7200   0.1674
  0.7400   0.1952
  0.7600   0.1941
  0.7800   0.1656
  0.8000   0.1152
  0.8200   0.0497
  0.8400  -0.0233
  0.8600  -0.0946
  0.8800  -0.1534
  0.9000  -0.1900
  0.9200  -0.1991
  0.9400  -0.1803
  0.9600  -0.1365
  0.9800  -0.0735
  1.0000  -0.0000
/\relax}\relax
\endpicture
}

\begin{figure}[h!]
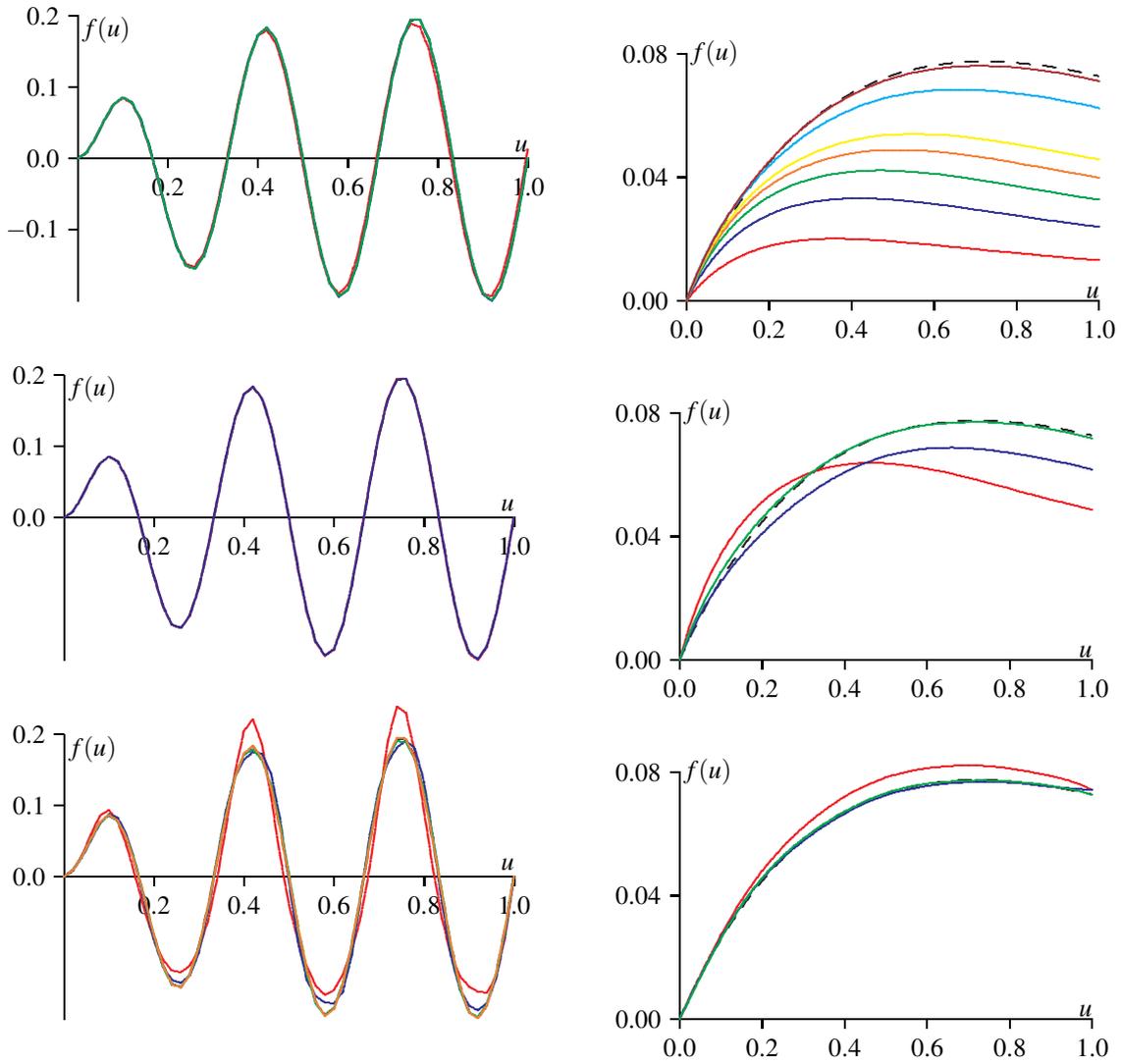

\hbox to \hsize{\hss\copy\figurefour\hss\copy\figureone\hss}
\vskip10pt
\hbox to \hsize{\hss\copy\figurefive\hss\copy\figuretwo\hss}
\vskip10pt
\hbox to \hsize{\hss\copy\figuresix\hss\copy\figurethree\hss}
\vskip10pt
\caption{\small {\bf Reconstructions of $f(u)$ from time trace (left) and final time  (right) data with Picard (top row), Anderson accelerated Picard (middle row), frozen Newton (bottom row); 
first (red), second (blue), third (green), fourth (orange) iterate versus exact $f$ (back dashed)}
}
\label{fig:f_W_1}
\end{figure}

It is well-known in undetermined coefficient problems that recovering a function of an independent variable from data in an orthogonal variable direction leads to severe ill-posedness whereas data given in a parallel direction often yields only mild ill-conditioning. In the case of the unknown depending on the state variable, here $u(x,t)$, one might assume that little difference would be found provided both data sets meet the range condition. However, previous work on such problems has shown that this is not always the case. In particular, for a parabolic type equation time trace data frequently leads to superior reconstructions over measurement data obtained in spatial manner
see, e.g., \cite{KaltenbacherRundell:2020b}. As we will see below, this is also the situation here. As a result, one is able to reconstruct more complex functions using time trace data for a given level of accuracy.
In selecting functions for testing we concentrated on highlighting this difference rather than functions that might correspond to a known physical reality. 
The chosen functions were
\begin{equation}\label{eqn:fact_titr}
f_{{\rm act}}(u) =  0.2\sin(6\pi u)(1-\exp(u))
\end{equation}
for time trace and
\begin{equation}\label{eqn:fact_fiti}
f_{{\rm act}}(u) =  0.1(1-\exp(-3u))\cos(0.7u)
\end{equation}
for final time data.

\begin{figure}[h!]
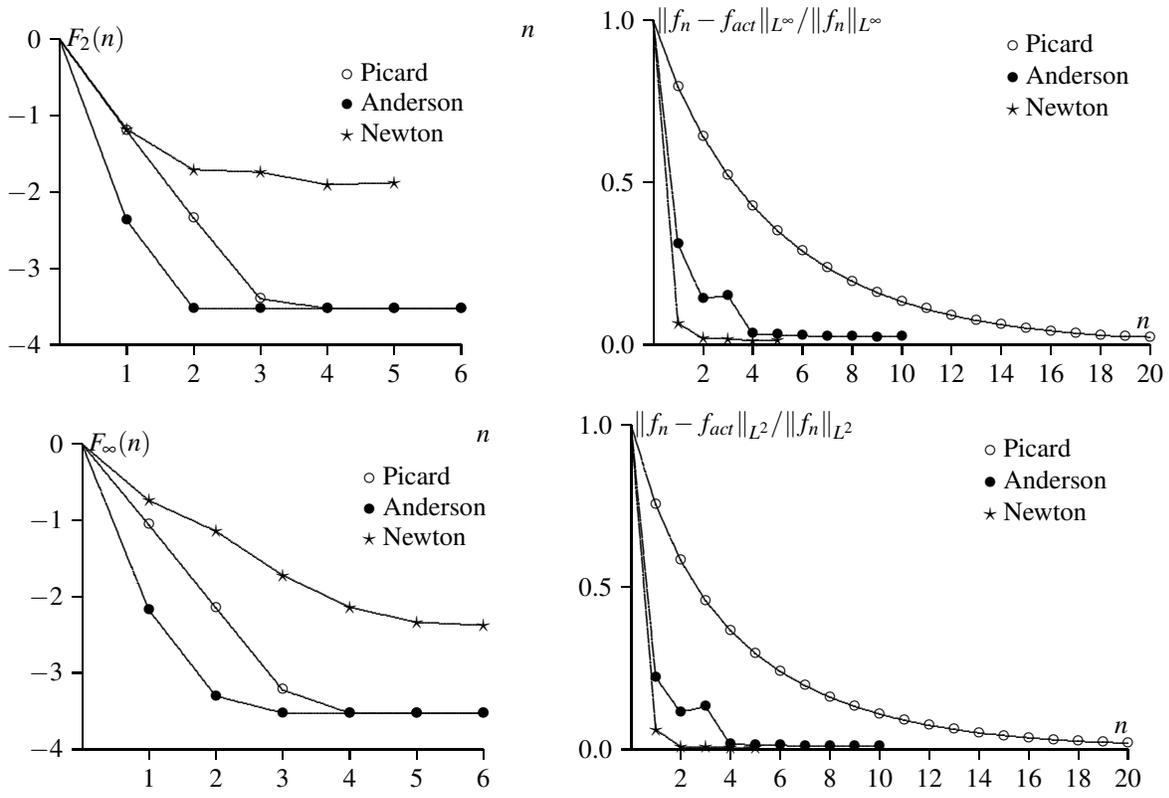

\hbox to \hsize{\hss\copy\figurethree\hss\copy\figurefour\hss}
\vskip10pt
\hbox to \hsize{\hss\copy\figurefive\hss\copy\figuresix\hss}
\vskip10pt
\caption{\small {\bf Relative error norms of the iterates for $f(u)$ from time trace (left) and final time (right) data.\\
$F_2(n) = \log_{10}\bigl(\|f_n-f_{act}\|_{L^2}/\|f_n\|_{L^2}\bigr)$ (top)
and $\,F_\infty(n) = \log_{10}\bigl(\|f_n-f_{act}\|_{L^\infty}/\|f_n\|_{L^\infty}\bigr)$ (bottom)
}}
\label{fig:f_W_2}
\end{figure}

In Figure \oldref{fig:f_W_1} we show reconstructions of the functions \eqref{eqn:fact_titr} and \eqref{eqn:fact_fiti}, while Figure \oldref{fig:f_W_2} displays the convergence history of the three different methods we use.

\font\smallsymbol = cmmi8
\newdimen\xfiglen \newdimen\yfiglen
\xfiglen=0.35 true in
\yfiglen=0.4 true in

\setbox\figurethree=\vbox{\hsize=\xfiglen
\beginpicture
\eightrm
\footnotesize
  \setcoordinatesystem units <\xfiglen,\yfiglen>  point at  0 -4
  \setplotarea x from 0 to 6, y from -4 to 0
  \axis bottom shiftedto y=-4 ticks short numbered from 1 to 6 by 1 /
  \axis left ticks short numbered from -4 to 0 by 1 /
%
\footnotesize
\put {$n$} [b] at 7 0.03
\put {$F_2(n)$} [l] at 0.1 0
\multiput {$\circ$} at 
    1.0000   -1.1952
    2.0000   -2.3372
    3.0000   -3.3979
    4.0000   -3.5229
    5.0000   -3.5229
    6.0000   -3.5229
/
\setlinear
\plot  
    0         0
    1.0000   -1.1952
    2.0000   -2.3372
    3.0000   -3.3979
    4.0000   -3.5229
    5.0000   -3.5229
    6.0000   -3.5229
/
\multiput {$\bullet$} at 
    1.0000   -2.3665
    2.0000   -3.5229
    3.0000   -3.5229
    4.0000   -3.5229
    5.0000   -3.5229
    6.0000   -3.5229
/
\setlinear
\plot  
    0.0000    0.0000
    1.0000   -2.3665
    2.0000   -3.5229
    3.0000   -3.5229
    4.0000   -3.5229
    5.0000   -3.5229
    6.0000   -3.5229
/
\multiput {$\star$} at
    1.0000   -1.1746
    2.0000   -1.7077
    3.0000   -1.7399
    4.0000   -1.9031
    5.0000   -1.8761
/
\setlinear
\plot  
  0.0000   0.0000
    1.0000   -1.1746
    2.0000   -1.7077
    3.0000   -1.7399
    4.0000   -1.9031
    5.0000   -1.8761
/
\footnotesize
\put {$\circ$ Picard} [lt] at 4.2 -0.3
\put {$\bullet$ Anderson} [lt] at 4.2 -0.7
\put {$\star$ Newton} [lt] at 4.2 -1.1
\endpicture
}
\xfiglen=0.13 true in
\yfiglen=1.7 true in
\setbox\figurefour=\vbox{\hsize=\xfiglen
\beginpicture
\eightrm
\footnotesize
  \setcoordinatesystem units <\xfiglen,\yfiglen>  point at  0 0
  \setplotarea x from 0 to 20, y from 0 to 1
  \axis bottom shiftedto y=0 ticks short numbered from 2 to 20 by 2 /
  \axis left ticks short numbered from 0 to 1 by 0.5 /
%
\footnotesize
\put {$n$} [rb] at 20 0.05
\put {$\|f_n-f_{act}\|_{L^\infty}/\|f_n\|_{L^\infty}$} [l] at 0.1 1
\multiput {$\circ$} at 
  1.0000   0.7952
  2.0000   0.6424
  3.0000   0.5234
  4.0000   0.4283
  5.0000   0.3514
  6.0000   0.2888
  7.0000   0.2376
  8.0000   0.1956
  9.0000   0.1612
 10.0000   0.1330
 11.0000   0.1098
 12.0000   0.0907
 13.0000   0.0751
 14.0000   0.0622
 15.0000   0.0516
 16.0000   0.0429
 17.0000   0.0357
 18.0000   0.0298
 19.0000   0.0249
 20.0000   0.0237
/
\setlinear
\plot  
  0.0000   1.0000
  1.0000   0.7952
  2.0000   0.6424
  3.0000   0.5234
  4.0000   0.4283
  5.0000   0.3514
  6.0000   0.2888
  7.0000   0.2376
  8.0000   0.1956
  9.0000   0.1612
 10.0000   0.1330
 11.0000   0.1098
 12.0000   0.0907
 13.0000   0.0751
 14.0000   0.0622
 15.0000   0.0516
 16.0000   0.0429
 17.0000   0.0357
 18.0000   0.0298
 19.0000   0.0249
 20.0000   0.0237
/
\multiput {$\bullet$} at 
  1.0000   0.3094
  2.0000   0.1415
  3.0000   0.1520
  4.0000   0.0336
  5.0000   0.0312
  6.0000   0.0288
  7.0000   0.0248
  8.0000   0.0242
  9.0000   0.0239
 10.0000   0.0245
/
\setlinear
\plot  
  0.0000   1.0000
  1.0000   0.3094
  2.0000   0.1415
  3.0000   0.1520
  4.0000   0.0336
  5.0000   0.0312
  6.0000   0.0288
  7.0000   0.0248
  8.0000   0.0242
  9.0000   0.0239
 10.0000   0.0245
/
\multiput {$\star$} at
  1.0000   0.0669
  2.0000   0.0196
  3.0000   0.0182
  4.0000   0.0125
  5.0000   0.0133
/
\setlinear
\plot  
  0.0000   1.0000
  1.0000   0.0669
  2.0000   0.0196
  3.0000   0.0182
  4.0000   0.0125
  5.0000   0.0133
/
\footnotesize
\put {$\circ$ Picard} [lb] at 14.2 0.9
\put {$\bullet$ Anderson} [lt] at 14.2 0.86
\put {$\star$ Newton} [lt] at 14.2 0.76
\endpicture
}

\xfiglen=0.35 true in
\yfiglen=0.4 true in
\setbox\figurefive=\vbox{\hsize=\xfiglen
\beginpicture
\eightrm
\footnotesize
  \setcoordinatesystem units <\xfiglen,\yfiglen>  point at  0 -4
  \setplotarea x from 0 to 6, y from -4 to 0
  \axis bottom shiftedto y=-4 ticks short numbered from 1 to 6 by 1 /
  \axis left ticks short numbered from -4 to 0 by 1 /
%
\footnotesize
\put {$n$} [b] at 6 0.03
\put {$F_\infty(n)$} [l] at 0.1 0
\multiput {$\circ$} at 
    1.0000   -1.0600
    2.0000   -2.1427
    3.0000   -3.2218
    4.0000   -3.5229
    5.0000   -3.5229
    6.0000   -3.5229
/
\setlinear
\plot  
    0         0
    1.0000   -1.0600
    2.0000   -2.1427
    3.0000   -3.2218
    4.0000   -3.5229
    5.0000   -3.5229
    6.0000   -3.5229
/
\multiput {$\bullet$} at 
    1.0000   -2.1739
    2.0000   -3.3010
    3.0000   -3.5229
    4.0000   -3.5229
    5.0000   -3.5229
    6.0000   -3.5229
/
\setlinear
\plot  
    0         0
    1.0000   -2.1739
    2.0000   -3.3010
    3.0000   -3.5229
    4.0000   -3.5229
    5.0000   -3.5229
    6.0000   -3.5229
/
\multiput {$\star$} at
    1.0000   -0.7409
    2.0000   -1.1415
    3.0000   -1.7190
    4.0000   -2.1367
    5.0000   -2.3372
    6.0000   -2.3768
/
\setlinear
\plot  
    0         0
    1.0000   -0.7409
    2.0000   -1.1415
    3.0000   -1.7190
    4.0000   -2.1367
    5.0000   -2.3372
    6.0000   -2.3768
/
\footnotesize
\put {$\circ$ Picard} [lt] at 4.2 -0.3
\put {$\bullet$ Anderson} [lt] at 4.2 -0.7
\put {$\star$ Newton} [lt] at 4.2 -1.1
\endpicture
}
\xfiglen=0.13 true in
\yfiglen=1.7 true in
\setbox\figuresix=\vbox{\hsize=\xfiglen
\beginpicture
\eightrm
\footnotesize
  \setcoordinatesystem units <\xfiglen,\yfiglen>  point at  0 0
  \setplotarea x from 0 to 20, y from 0 to 1
  \axis bottom shiftedto y=0 ticks short numbered from 2 to 20 by 2 /
  \axis left ticks short numbered from 0 to 1 by 0.5 /
%
\footnotesize
\put {$n$} [rb] at 20 0.05
\put {$\|f_n-f_{act}\|_{L^2}/\|f_n\|_{L^2}$} [l] at 0.1 1
\multiput {$\circ$} at 
  1.0000   0.7542
  2.0000   0.5824
  3.0000   0.4586
  4.0000   0.3665
  5.0000   0.2959
  6.0000   0.2406
  7.0000   0.1966
  8.0000   0.1611
  9.0000   0.1324
 10.0000   0.1090
 11.0000   0.0898
 12.0000   0.0742
 13.0000   0.0615
 14.0000   0.0510
 15.0000   0.0425
 16.0000   0.0356
 17.0000   0.0299
 18.0000   0.0253
 19.0000   0.0216
 20.0000   0.0186
/
\setlinear
\plot  
  0.0000   1.0000
  1.0000   0.7542
  2.0000   0.5824
  3.0000   0.4586
  4.0000   0.3665
  5.0000   0.2959
  6.0000   0.2406
  7.0000   0.1966
  8.0000   0.1611
  9.0000   0.1324
 10.0000   0.1090
 11.0000   0.0898
 12.0000   0.0742
 13.0000   0.0615
 14.0000   0.0510
 15.0000   0.0425
 16.0000   0.0356
 17.0000   0.0299
 18.0000   0.0253
 19.0000   0.0216
 20.0000   0.0186
/
\multiput {$\bullet$} at 
  1.0000   0.2224
  2.0000   0.1148
  3.0000   0.1324
  4.0000   0.0159
  5.0000   0.0139
  6.0000   0.0124
  7.0000   0.0102
  8.0000   0.0098
  9.0000   0.0097
 10.0000   0.0099
/
\setlinear
\plot  
  0.0000   1.0000
  1.0000   0.2224
  2.0000   0.1148
  3.0000   0.1324
  4.0000   0.0159
  5.0000   0.0139
  6.0000   0.0124
  7.0000   0.0102
  8.0000   0.0098
  9.0000   0.0097
 10.0000   0.0099
/
\multiput {$\star$} at
  1.0000   0.0608
  2.0000   0.0088
  3.0000   0.0086
  4.0000   0.0057
  5.0000   0.0061
/
\setlinear
\plot  
  0.0000   1.0000
  1.0000   0.0608
  2.0000   0.0088
  3.0000   0.0086
  4.0000   0.0057
  5.0000   0.0061
/
\footnotesize
\put {$\circ$ Picard} [lb] at 14.2 0.9
\put {$\bullet$ Anderson} [lt] at 14.2 0.86
\put {$\star$ Newton} [lt] at 14.2 0.76
\endpicture
}
%

\begin{figure}[h!]
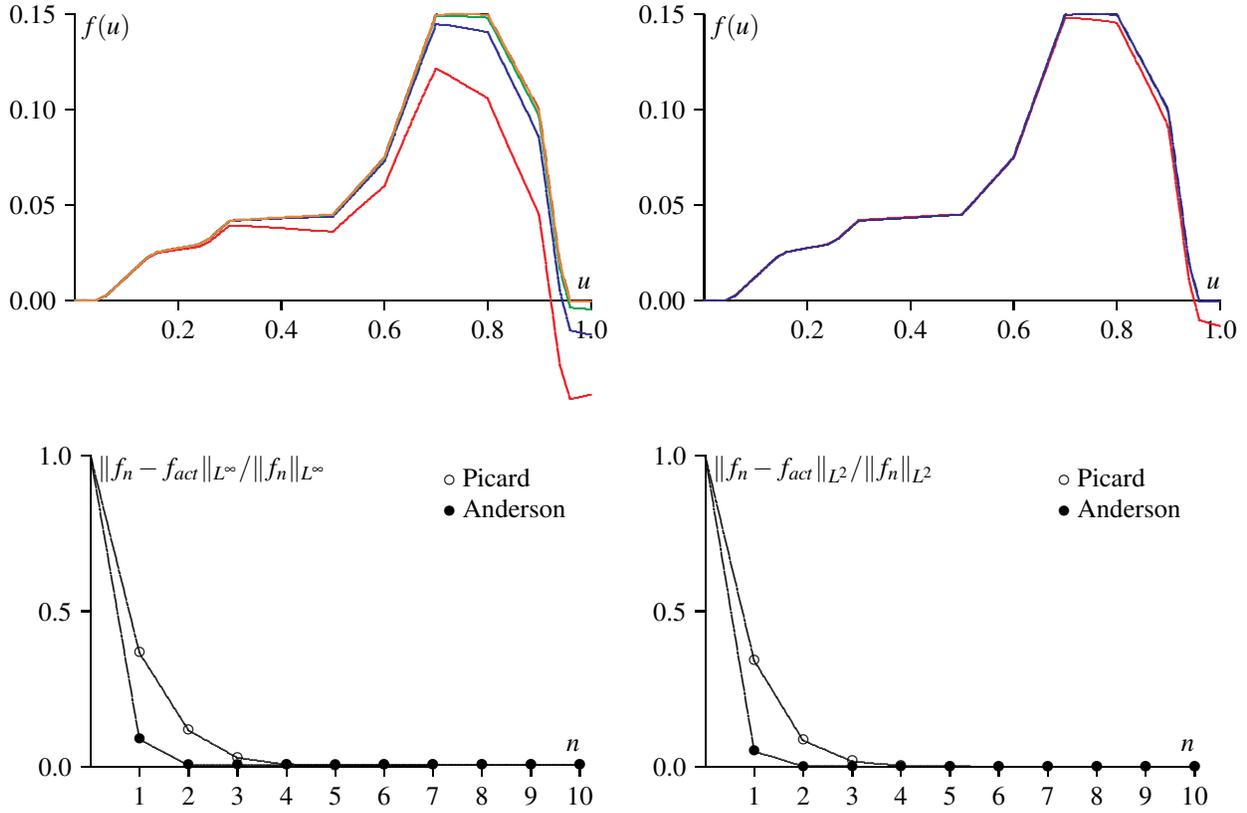

\hbox to \hsize{\hss\copy\figureseven\hss\copy\figureeight\hss}
\vskip40pt
\hbox to \hsize{\hss\copy\figurenine\hss\copy\figureten\hss}
\vskip10pt
\caption{\small {\bf Reconstructions of $f(u)$ for piecewise linear example from time trace data with Picard (top left) and Anderson accelerated Picard (top right); 
first (red), second (blue), third (green), fourth (orange), fifth (yellow), sixth (cyan), seventh (maroon) iterate versus exact $f$ (back dashed);
relative error norms of the iterates for $f(u)$ in $L^\infty$ (bottom left) and $L^2$ (bottom right) norm}}
\label{fig:f_W_3}
\end{figure}

The two examples shown above represent smooth functions and a reasonable question is what is the ability to reconstruction a function lying at the boundary of the regularity theorems in section \oldref{sec:forward}, namely a function that is just Lipschitz continuous.
In Figure \oldref{fig:f_W_3} we therefore show the results of the fixed point scheme and its Anderson accelerated version using time trace data. For doing this with Newton, the optimal basis would probably have to be piecewise linear, as for example in \cite{KaltenbacherRundell:2021d}. 
To avoid a mismatch between the actual location of the non-differentiability points and the nodes of the basis functions requires a large number of basis functions in general, unless these points are a priori known. Thus without this knowledge, the computational complexity of calculating the Jacobian would make (frozen) Newton much slower than the accelerated fixed point schemes.  

\xfiglen=2.7 true in
\yfiglen=10 true in
\setbox\figureseven=\vbox{\hsize=\xfiglen
\beginpicture
\eightrm
\footnotesize
  \setcoordinatesystem units <\xfiglen,\yfiglen>  
  \setplotarea x from 0 to 1, y from 0 to 0.15
  \axis bottom shiftedto y=0 ticks short numbered from 0.2 to 1 by 0.2 /
  \axis left ticks short numbered from 0 to 0.15 by 0.05 /
\footnotesize
\put {$u$} [rb] at 1 0.006
\put {$f(u)$} [lt] at 0.013 0.15
\small
\linethickness=0.6pt
\setplotsymbol ({\sevenrm .})
\setquadratic
\setlinear
\setdashes
\Black{\relax  
\plot   
  0.0000   0.0000
  0.0200   0.0000
  0.0400   0.0000
  0.0600   0.0025
  0.0800   0.0075
  0.1000   0.0125
  0.1200   0.0175
  0.1400   0.0225
  0.1600   0.0255
  0.1800   0.0265
  0.2000   0.0275
  0.2200   0.0285
  0.2400   0.0295
  0.2600   0.0324
  0.2800   0.0372
  0.3000   0.0420
  0.3200   0.0423
  0.3400   0.0426
  0.3600   0.0429
  0.3800   0.0432
  0.4000   0.0435
  0.4200   0.0438
  0.4400   0.0441
  0.4600   0.0444
  0.4800   0.0447
  0.5000   0.0450
  0.5200   0.0510
  0.5400   0.0570
  0.5600   0.0630
  0.5800   0.0690
  0.6000   0.0750
  0.6200   0.0900
  0.6400   0.1050
  0.6600   0.1200
  0.6800   0.1350
  0.7000   0.1500
  0.7200   0.1500
  0.7400   0.1500
  0.7600   0.1500
  0.7800   0.1500
  0.8000   0.1500
  0.8200   0.1400
  0.8400   0.1300
  0.8600   0.1200
  0.8800   0.1100
  0.9000   0.1000
  0.9200   0.0600
  0.9400   0.0200
  0.9600   0.0000
  0.9800   0.0000
  1.0000   0.0000
/\relax}\relax
\setsolid
\Red{\relax 
\plot
  0.0000   0.0000
  0.0200   0.0000
  0.0400   0.0000
  0.0600   0.0025
  0.0800   0.0075
  0.1000   0.0124
  0.1200   0.0173
  0.1400   0.0222
  0.1600   0.0251
  0.1800   0.0259
  0.2000   0.0267
  0.2200   0.0274
  0.2400   0.0281
  0.2600   0.0307
  0.2800   0.0351
  0.3000   0.0394
  0.3200   0.0392
  0.3400   0.0390
  0.3600   0.0387
  0.3800   0.0384
  0.4000   0.0381
  0.4200   0.0377
  0.4400   0.0373
  0.4600   0.0369
  0.4800   0.0365
  0.5000   0.0361
  0.5200   0.0411
  0.5400   0.0460
  0.5600   0.0508
  0.5800   0.0554
  0.6000   0.0601
  0.6200   0.0731
  0.6400   0.0858
  0.6600   0.0983
  0.6800   0.1103
  0.7000   0.1217
  0.7200   0.1188
  0.7400   0.1156
  0.7600   0.1124
  0.7800   0.1093
  0.8000   0.1060
  0.8200   0.0934
  0.8400   0.0810
  0.8600   0.0689
  0.8800   0.0571
  0.9000   0.0450
  0.9200   0.0054
  0.9400  -0.0335
  0.9600  -0.0517
  0.9800  -0.0504
  1.0000  -0.0491
/\relax}\relax
\Blue{\relax 
\plot
  0.0000   0.0000
  0.0200   0.0000
  0.0400   0.0000
  0.0600   0.0025
  0.0800   0.0075
  0.1000   0.0125
  0.1200   0.0175
  0.1400   0.0225
  0.1600   0.0255
  0.1800   0.0265
  0.2000   0.0275
  0.2200   0.0285
  0.2400   0.0295
  0.2600   0.0323
  0.2800   0.0371
  0.3000   0.0418
  0.3200   0.0421
  0.3400   0.0424
  0.3600   0.0426
  0.3800   0.0429
  0.4000   0.0431
  0.4200   0.0433
  0.4400   0.0435
  0.4600   0.0437
  0.4800   0.0438
  0.5000   0.0441
  0.5200   0.0498
  0.5400   0.0556
  0.5600   0.0613
  0.5800   0.0671
  0.6000   0.0730
  0.6200   0.0874
  0.6400   0.1019
  0.6600   0.1164
  0.6800   0.1308
  0.7000   0.1447
  0.7200   0.1444
  0.7400   0.1436
  0.7600   0.1427
  0.7800   0.1418
  0.8000   0.1406
  0.8200   0.1299
  0.8400   0.1189
  0.8600   0.1079
  0.8800   0.0969
  0.9000   0.0853
  0.9200   0.0456
  0.9400   0.0052
  0.9600  -0.0155
  0.9800  -0.0167
  1.0000  -0.0177
/\relax}\relax
\Green{\relax 
\plot
  0.0000   0.0000
  0.0200   0.0000
  0.0400   0.0000
  0.0600   0.0025
  0.0800   0.0075
  0.1000   0.0125
  0.1200   0.0175
  0.1400   0.0225
  0.1600   0.0255
  0.1800   0.0265
  0.2000   0.0275
  0.2200   0.0285
  0.2400   0.0295
  0.2600   0.0324
  0.2800   0.0372
  0.3000   0.0419
  0.3200   0.0423
  0.3400   0.0426
  0.3600   0.0429
  0.3800   0.0432
  0.4000   0.0435
  0.4200   0.0438
  0.4400   0.0441
  0.4600   0.0443
  0.4800   0.0446
  0.5000   0.0450
  0.5200   0.0509
  0.5400   0.0569
  0.5600   0.0628
  0.5800   0.0688
  0.6000   0.0749
  0.6200   0.0897
  0.6400   0.1046
  0.6600   0.1195
  0.6800   0.1344
  0.7000   0.1490
  0.7200   0.1491
  0.7400   0.1490
  0.7600   0.1488
  0.7800   0.1486
  0.8000   0.1482
  0.8200   0.1381
  0.8400   0.1279
  0.8600   0.1176
  0.8800   0.1073
  0.9000   0.0964
  0.9200   0.0569
  0.9400   0.0167
  0.9600  -0.0036
  0.9800  -0.0040
  1.0000  -0.0044
/\relax}\relax
\Orange{\relax 
\plot
  0.0000   0.0000
  0.0200   0.0000
  0.0400   0.0000
  0.0600   0.0025
  0.0800   0.0075
  0.1000   0.0125
  0.1200   0.0175
  0.1400   0.0225
  0.1600   0.0255
  0.1800   0.0265
  0.2000   0.0275
  0.2200   0.0285
  0.2400   0.0295
  0.2600   0.0324
  0.2800   0.0372
  0.3000   0.0419
  0.3200   0.0423
  0.3400   0.0426
  0.3600   0.0429
  0.3800   0.0432
  0.4000   0.0435
  0.4200   0.0438
  0.4400   0.0441
  0.4600   0.0444
  0.4800   0.0447
  0.5000   0.0451
  0.5200   0.0510
  0.5400   0.0570
  0.5600   0.0630
  0.5800   0.0690
  0.6000   0.0752
  0.6200   0.0900
  0.6400   0.1050
  0.6600   0.1200
  0.6800   0.1350
  0.7000   0.1497
  0.7200   0.1500
  0.7400   0.1500
  0.7600   0.1500
  0.7800   0.1500
  0.8000   0.1498
  0.8200   0.1400
  0.8400   0.1300
  0.8600   0.1199
  0.8800   0.1099
  0.9000   0.0993
  0.9200   0.0599
  0.9400   0.0199
  0.9600  -0.0001
  0.9800  -0.0001
  1.0000  -0.0002
/\relax}\relax
\endpicture
}

\setbox\figureeight=\vbox{\hsize=\xfiglen
\beginpicture
\eightrm
\footnotesize
  \setcoordinatesystem units <\xfiglen,\yfiglen>  
  \setplotarea x from 0 to 1, y from 0 to 0.15
  \axis bottom shiftedto y=0 ticks short numbered from 0.2 to 1 by 0.2 /
  \axis left ticks short numbered from 0 to 0.15 by 0.05 /
\footnotesize
\put {$u$} [rb] at 1 0.006
\put {$f(u)$} [lt] at 0.013 0.15
\small
\linethickness=0.6pt
\setplotsymbol ({\sevenrm .})
\setquadratic
\setlinear
\setdashes
\Black{\relax  
\plot   
  0.0000   0.0000
  0.0200   0.0000
  0.0400   0.0000
  0.0600   0.0025
  0.0800   0.0075
  0.1000   0.0125
  0.1200   0.0175
  0.1400   0.0225
  0.1600   0.0255
  0.1800   0.0265
  0.2000   0.0275
  0.2200   0.0285
  0.2400   0.0295
  0.2600   0.0324
  0.2800   0.0372
  0.3000   0.0420
  0.3200   0.0423
  0.3400   0.0426
  0.3600   0.0429
  0.3800   0.0432
  0.4000   0.0435
  0.4200   0.0438
  0.4400   0.0441
  0.4600   0.0444
  0.4800   0.0447
  0.5000   0.0450
  0.5200   0.0510
  0.5400   0.0570
  0.5600   0.0630
  0.5800   0.0690
  0.6000   0.0750
  0.6200   0.0900
  0.6400   0.1050
  0.6600   0.1200
  0.6800   0.1350
  0.7000   0.1500
  0.7200   0.1500
  0.7400   0.1500
  0.7600   0.1500
  0.7800   0.1500
  0.8000   0.1500
  0.8200   0.1400
  0.8400   0.1300
  0.8600   0.1200
  0.8800   0.1100
  0.9000   0.1000
  0.9200   0.0600
  0.9400   0.0200
  0.9600   0.0000
  0.9800   0.0000
  1.0000   0.0000
/\relax}\relax
\setsolid
\Red{\relax 
\plot
  0.0000   0.0000
  0.0200   0.0000
  0.0400   0.0000
  0.0600   0.0025
  0.0800   0.0075
  0.1000   0.0125
  0.1200   0.0175
  0.1400   0.0225
  0.1600   0.0256
  0.1800   0.0266
  0.2000   0.0276
  0.2200   0.0286
  0.2400   0.0297
  0.2600   0.0326
  0.2800   0.0374
  0.3000   0.0422
  0.3200   0.0426
  0.3400   0.0429
  0.3600   0.0432
  0.3800   0.0435
  0.4000   0.0438
  0.4200   0.0441
  0.4400   0.0444
  0.4600   0.0447
  0.4800   0.0449
  0.5000   0.0453
  0.5200   0.0511
  0.5400   0.0570
  0.5600   0.0629
  0.5800   0.0687
  0.6000   0.0748
  0.6200   0.0894
  0.6400   0.1042
  0.6600   0.1190
  0.6800   0.1337
  0.7000   0.1480
  0.7200   0.1480
  0.7400   0.1476
  0.7600   0.1471
  0.7800   0.1465
  0.8000   0.1456
  0.8200   0.1352
  0.8400   0.1244
  0.8600   0.1136
  0.8800   0.1027
  0.9000   0.0911
  0.9200   0.0514
  0.9400   0.0108
  0.9600  -0.0102
  0.9800  -0.0118
  1.0000  -0.0131
/\relax}\relax
\Blue{\relax 
\plot
  0.0000   0.0000
  0.0200   0.0000
  0.0400   0.0000
  0.0600   0.0025
  0.0800   0.0075
  0.1000   0.0125
  0.1200   0.0175
  0.1400   0.0225
  0.1600   0.0255
  0.1800   0.0265
  0.2000   0.0275
  0.2200   0.0285
  0.2400   0.0295
  0.2600   0.0324
  0.2800   0.0372
  0.3000   0.0419
  0.3200   0.0422
  0.3400   0.0425
  0.3600   0.0428
  0.3800   0.0431
  0.4000   0.0434
  0.4200   0.0437
  0.4400   0.0441
  0.4600   0.0444
  0.4800   0.0447
  0.5000   0.0451
  0.5200   0.0510
  0.5400   0.0570
  0.5600   0.0630
  0.5800   0.0690
  0.6000   0.0752
  0.6200   0.0901
  0.6400   0.1051
  0.6600   0.1201
  0.6800   0.1351
  0.7000   0.1498
  0.7200   0.1500
  0.7400   0.1500
  0.7600   0.1500
  0.7800   0.1500
  0.8000   0.1498
  0.8200   0.1400
  0.8400   0.1300
  0.8600   0.1199
  0.8800   0.1099
  0.9000   0.0992
  0.9200   0.0597
  0.9400   0.0197
  0.9600  -0.0003
  0.9800  -0.0004
  1.0000  -0.0005
/\relax}\relax
\endpicture
}

\xfiglen=0.32 true in
\yfiglen=1.63 true in
\setbox\figurenine=\vbox{\hsize=\xfiglen
\beginpicture
\eightrm
\footnotesize
  \setcoordinatesystem units <0.8\xfiglen,\yfiglen> 
  \setplotarea x from 0 to 7, y from 0 to 1
  \axis bottom shiftedto y=0 ticks short numbered from 1 to 10 by 1 /
  \axis left ticks short numbered from 0 to 1 by 0.5 /
%
\footnotesize
\put {$n$} [rb] at 10 0.05
\put {$\|f_n-f_{act}\|_{L^\infty}/\|f_n\|_{L^\infty}$} [lt] at 0.15 1
\multiput {$\circ$} at 
  1.0000   0.3669
  2.0000   0.1178
  3.0000   0.0293
  4.0000   0.0077
  5.0000   0.0047
  6.0000   0.0051
  7.0000   0.0052
  8.0000   0.0051
  9.0000   0.0051
 10.0000   0.0050
/
\setlinear
\plot  
  0.0000   1.0000
  1.0000   0.3669
  2.0000   0.1178
  3.0000   0.0293
  4.0000   0.0077
  5.0000   0.0047
  6.0000   0.0051
  7.0000   0.0052
  8.0000   0.0051
  9.0000   0.0051
 10.0000   0.0050
/
\multiput {$\bullet$} at 
  1.0000   0.0877
  2.0000   0.0053
  3.0000   0.0053
  4.0000   0.0053
  5.0000   0.0052
  6.0000   0.0051
  7.0000   0.0050
  8.0000   0.0049
  9.0000   0.0049
 10.0000   0.0049
/
\setlinear
\plot  
  0.0000   1.0000
  1.0000   0.0877
  2.0000   0.0053
  3.0000   0.0053
  4.0000   0.0053
  5.0000   0.0052
  6.0000   0.0051
  7.0000   0.0050
  8.0000   0.0049
  9.0000   0.0049
 10.0000   0.0049
/
\footnotesize
\put {$\circ$ Picard} [lb] at 7.2 0.9
\put {$\bullet$ Anderson} [lt] at 7.2 0.86
\endpicture
}
\setbox\figureten=\vbox{\hsize=\xfiglen
\beginpicture
\eightrm
\footnotesize
  \setcoordinatesystem units <0.8\xfiglen,\yfiglen> 
  \setplotarea x from 0 to 10, y from 0 to 1
  \axis bottom shiftedto y=0 ticks short numbered from 1 to 10 by 1 /
  \axis left ticks short numbered from 0 to 1 by 0.5 /
%
\footnotesize
\put {$n$} [rb] at 10 0.05
\put {$\|f_n-f_{act}\|_{L^2}/\|f_n\|_{L^2}$} [lt] at 0.15 1
\multiput {$\circ$} at 
  1.0000   0.3418
  2.0000   0.0849
  3.0000   0.0178
  4.0000   0.0034
  5.0000   0.0011
  6.0000   0.0010
  7.0000   0.0010
  8.0000   0.0010
  9.0000   0.0010
 10.0000   0.0009
/
\setlinear
\plot  
  0.0000   1.0000
  1.0000   0.3418
  2.0000   0.0849
  3.0000   0.0178
  4.0000   0.0034
  5.0000   0.0011
  6.0000   0.0010
  7.0000   0.0010
  8.0000   0.0010
  9.0000   0.0010
 10.0000   0.0009
/
\multiput {$\bullet$} at 
  1.0000   0.0501
  2.0000   0.0017
  3.0000   0.0010
  4.0000   0.0010
  5.0000   0.0010
  6.0000   0.0010
  7.0000   0.0009
  8.0000   0.0009
  9.0000   0.0009
 10.0000   0.0009
/
\setlinear
\plot  
  0.0000   1.0000
  1.0000   0.0501
  2.0000   0.0017
  3.0000   0.0010
  4.0000   0.0010
  5.0000   0.0010
  6.0000   0.0010
  7.0000   0.0009
  8.0000   0.0009
  9.0000   0.0009
 10.0000   0.0009
/
\footnotesize
\put {$\circ$ Picard} [lb] at 7.2 0.9
\put {$\bullet$ Anderson} [lt] at 7.2 0.86
\endpicture
}
%

We also did reconstructions from time trace data in the fractional case \eqref{Westervelt_f_int_frac2}, where as mentioned in section \oldref{sec:reconstruction_schemes} 
the time trace data iteration scheme \eqref{iteration_titr} remains the same as in case $\alpha=1$, with the dependence on $\alpha$ showing up only in the {\sc pde} to be solved.
Provided $\alpha$ was not close to zero, the differences in the reconstructions were very little.

\subsection*{Appendix (Proof of Lemma \ref{lem:enest_lin})}

Multiplying \eqref{Westervelt_u_lin} with $\ppsi_{tt}$ and integrating over $\Omega\times(0,t)$ yields 
\begin{equation}\label{eqn:enid_1-d}
\begin{aligned}
&\int_0^t \|\sqrt{1-\sigma(\tau)}\ppsi_{tt}(\tau)\|_{L^2(\Omega)}^2\, d\tau 
+ \frac{b}{2}\|\nabla \ppsi_t(t)\|_{L^2(\Omega)}^2
\\
&= \frac{b}{2}\|\nabla \ppsi_1\|_{L^2(\Omega)}^2 
+ c^2 \int_0^t \|\nabla \ppsi_t(\tau)\|_{L^2(\Omega)}^2\, d\tau 
- c^2 \langle \nabla \ppsi_t(t), \nabla \ppsi(t) \rangle_{L^2(\Omega)}
\\&\qquad
+\int_0^t\langle \tilde{r}(\tau)+\eta(\tau)\ppsi_t(\tau),\ppsi_{tt}(\tau)\rangle_{L^2(\Omega)}\, d\tau\\
&\leq 
\frac{b}{2}\|\nabla \ppsi_1\|_{L^2(\Omega)}^2  
+ c^2 \int_0^t \|\nabla \ppsi_t(\tau)\|_{L^2(\Omega)}^2\, d\tau
+\frac{b}{4} \|\nabla \ppsi_t(t)\|_{L^2(\Omega)}^2
+ \frac{c^4T}{b}\int_0^t\|\nabla \ppsi_t(\tau)\, d\tau\|_{L^2(\Omega)}^2 \\
&\qquad 
+\frac12 \int_0^t \|\sqrt{1-\sigma(\tau)}\ppsi_{tt}(\tau)\|_{L^2(\Omega)}^2\, d\tau 
+ \frac{1}{2(1-\overline{\sigma})}\int_0^t\|\tilde{r}(\tau)+\eta(\tau)\ppsi_t(\tau)\|_{L^2(\Omega)}^2\, d\tau,
\end{aligned}
\end{equation}
where we have used $\sigma\leq\overline{\sigma}<1$, $\ppsi(0)=0$, 
$\|\nabla\ppsi(t)\|_{L^2(\Omega)}^2 = \|\int_0^t\nabla\ppsi_t(\tau)\, d\tau\|_{L^2(\Omega)}^2
\leq t\int_0^t\|\nabla\ppsi_t(\tau)\|_{L^2(\Omega)}^2\, d\tau$ and Young's inequality.
From this we extract 
\[
\begin{aligned}
&\frac{1-\overline{\sigma}}{2}\int_0^t\|\ppsi_{tt}(\tau)\|_{L^2(\Omega)}^2\, d\tau
+\frac{b}{4}\|\nabla \ppsi_t(t)\|_{L^2(\Omega)}^2\\
&\leq \frac{b}{2}\|\nabla \ppsi_1\|_{L^2(\Omega)}^2 
+ \frac{1}{2(1-\overline{\sigma})}\int_0^t\|\tilde{r}(\tau)+\eta(\tau)\ppsi_t(\tau)\|_{L^2(\Omega)}^2\, d\tau\\
&\qquad+ \Bigl(c^2+\frac{c^4T}{b}\Bigr) \int_0^t \|\nabla \ppsi_t(\tau)\|_{L^2(\Omega)}^2\, d\tau\,.
\end{aligned}
\]
To estimate the term containing $\eta$ we make use of continuity of the embeddings $H^1(\Omega)\to L^\infty(\Omega)$ in one space dimension or $H^1(\Omega)\to L^\infty(\Omega)$ in two and three space dimensions, respectively, as well as the Poincar\'{e}-Friedrichs inequality 
\[
\|v\|_{L^p(\Omega)}\leq C_{H^1,L^p}^\Omega  \|\nabla v\|_{L^2(\Omega)} \mbox{ for all }v\in H_0^1(\Omega)\,,
\mbox{ with } \begin{cases}
p=\infty \ \mbox{ if }d=1\\
p\leq 6 \ \mbox{ if }d\in\{2,3\}\,.
\end{cases}
\]
to conclude
$\int_0^t\|\eta(\tau)\ppsi_t(\tau)\|_{L^2(\Omega)}^2\, d\tau
\leq C(\eta) \|\nabla\ppsi_t\|_{L^\infty(0,t;L^2(\Omega)}^2$
where we assume  $C(\eta)$ as in \eqref{eqn:etasmall} to be smaller than $\frac{b(1-\overline{\sigma})}{4}$.
Gronwall's inequality for 
\[
\xi(t)=\max\left\{\Bigl(\frac{b}{4}-\frac{C(\eta)}{1-\bar{\sigma}}\Bigr)\|\nabla \ppsi_t\|_{L^\infty(0,t;L^2(\Omega)}^2\, \ \frac{1-\bar{\sigma}}{2}\|\ppsi_{tt}\|_{L^2(0,t;L^2(\Omega)}^2\right\}
\]
satisfying 
\[
\xi(t) \leq \frac{b}{2}\|\nabla \ppsi_1\|_{L^2(\Omega)}^2 
+ \frac{1}{1-\overline{\sigma}}\|\tilde{r}\|_{L^2(0,T;L^2(\Omega))}^2
+ \Bigl(c^2+\frac{c^4T}{b}\Bigr) \frac{4(1-\overline{\sigma})}{b(1-\overline{\sigma})-4C(\eta)}
\int_0^t \xi(\tau)\, d\tau
\]
yields \eqref{eqn:enest_mult_utt}.
  
\medskip

To prove  \eqref{eqn:enest_mult_Autt}, we multiply \eqref{Westervelt_u_lin} with $\mathcal{A}\ppsi_{tt}$ and integrate over $\Omega\times(0,t)$, which similarly to \eqref{eqn:enid_1-d}, but taking into account
\[
\int_\Omega (1-\sigma) \ppsi_{tt} \mathcal{A} \ppsi_{tt}\, dx 
=\int_\Omega \Bigl((1-\sigma) |\nabla \ppsi_{tt}|^2 - \ppsi_{tt}\nabla\sigma\cdot\nabla \ppsi_{tt}\, dx 
\]
yields 
\begin{equation}\label{eqn:enid}
\begin{aligned}
\int_0^t \|\sqrt{1-\sigma(\tau)}&\nabla \ppsi_{tt}(\tau)\|_{L^2(\Omega)}^2\, d\tau 
+ \frac{b}{2}\|\mathcal{A} \ppsi_t(t)\|_{L^2(\Omega)}^2
\\
&= \frac{b}{2}\|\mathcal{A} \ppsi_1\|_{L^2(\Omega)}^2 
+ c^2 \int_0^t \|\mathcal{A} \ppsi_t(\tau)\|_{L^2(\Omega)}^2\, d\tau 
- c^2 \langle \mathcal{A} \ppsi_t(t), \nabla \ppsi(t) \rangle_{L^2(\Omega)}\\
&\qquad+\int_0^t\langle \ppsi_{tt}\nabla\sigma + \nabla (\tilde{r}+\eta\ppsi_t)(\tau),\nabla \ppsi_{tt}(\tau)\rangle_{L^2(\Omega)}\, d\tau\\
&\leq 
\frac{b}{2}\|\mathcal{A} \ppsi_1\|_{L^2(\Omega)}^2  
+ c^2 \int_0^t \|\mathcal{A} \ppsi_t(\tau)\|_{L^2(\Omega)}^2\, d\tau\\
&\qquad+\frac{1-\overline{\sigma}}{2} \int_0^t \|\nabla \ppsi_{tt}(\tau)\|_{L^2(\Omega)}^2\, d\tau 
+\frac{b}{4} \|\mathcal{A} \ppsi_t(t)\|_{L^2(\Omega)}^2 \\
&\qquad+ \frac{c^4}{b}\int_0^t\|\mathcal{A} \ppsi_t(\tau)\, d\tau\|_{L^2(\Omega)}^2 
+ \frac{1}{2(1-\overline{\sigma})}\int_0^t\|\nabla(\tilde{r}+\eta\ppsi_t)(\tau)\|_{L^2(\Omega)}^2\, d\tau\\
&\qquad+ \|\nabla\sigma\|_{L^\infty(0,T;L^3(\Omega))} \int_0^t \|\ppsi_{tt}(\tau)\|_{L^6(\Omega)}\|\nabla \ppsi_{tt}(\tau)\|_{L^2(\Omega)}\, d\tau.
\end{aligned}
\end{equation}
Using continuity of the embeddings $H_0^1(\Omega)\to L^6(\Omega)$ and $H^2(\Omega)\to L^\infty(\Omega)$ in 
\[
\begin{aligned}
&\int_0^t\|\nabla(\eta\ppsi_t)(\tau)\|_{L^2(\Omega)}^2 \, d\tau
\leq C(\eta)\|\mathcal{A}\ppsi_t\|_{L^\infty(0,t;L^2(\Omega)}^2\\
&\|\nabla\sigma\|_{L^\infty(0,T;L^3(\Omega))} \int_0^t \|\ppsi_{tt}(\tau)\|_{L^6(\Omega)}\|\nabla \ppsi_{tt}(\tau)\|_{L^2(\Omega)}\, d\tau \leq C(\sigma)\int_0^t \|\nabla \ppsi_{tt}(\tau)\|_{L^2(\Omega)}^2\, d\tau
\end{aligned}
\]
and imposing the smallness conditions \eqref{eqn:etasmall_mult_Autt}, \eqref{eqn:smallnablasigma} on $\eta$, $\nabla\eta$, and $\nabla\sigma$
analogously to above we obtain \eqref{eqn:enest_mult_Autt}.
\goodbreak
\section*{Acknowledgment}

\noindent
The work of the first author was supported by the Austrian Science Fund {\sc fwf}
under the grant P30054.

\noindent
The work of the second author was supported
in part by the
National Science Foundation through award {\sc dms}-2111020.

\noindent

\end{document}